\documentclass{article}
\usepackage[english]{babel}
\usepackage[latin1]{inputenc}
\usepackage[T1]{fontenc}
\usepackage{amsmath, amsfonts, amssymb, amsthm}
\newtheorem{theo}{Theorem}
\newtheorem{lemm}{Lemma}
\newtheorem{prop}{Proposition}
\newtheorem{rema}{Remark}
\newtheorem{coro}{Corollary}
\newtheorem{defi}{Definition}
\newcommand \e {\varepsilon}
\newcommand \ind {\mathbb{I}}

\begin{document}
\begin{center}
\textbf{QUENCHED CONVERGENCE OF A SEQUENCE OF SUPERPROCESSES IN
$\mathbb{R}^d$ AMONG POISSONIAN OBSTACLES}
\end{center}

\begin{center}
\textbf{Amandine V\'eber}
\end{center}

\begin{center}
D\'epartement de Math\'ematiques, Universit\'e Paris-Sud \\
91405 Orsay Cedex, France\\
\verb"amandine.veber@math.u-psud.fr"
\end{center}
\begin{abstract}
We prove a convergence theorem for a sequence of super-Brownian
motions moving among hard Poissonian obstacles, when the intensity
of the obstacles grows to infinity but their diameters shrink to
zero in an appropriate manner. The superprocesses are shown to
converge in probability for the law $\mathbf{P}$ of the obstacles,
and $\mathbf{P}$-almost surely for a subsequence, towards a
superprocess with underlying spatial motion given by Brownian motion
and (inhomogeneous) branching mechanism $\psi(u,x)$ of the form
$\psi(u,x)= u^2+ \kappa(x)u$, where $\kappa(x)$ depends on the
density of the obstacles. This work draws on similar questions for a
single Brownian motion. In the course of the proof, we establish
precise estimates for integrals of functions over the Wiener
sausage, which are of independent interest.


\end{abstract}

\section{Introduction}
\subsection{Superprocesses in random media}
The purpose of this article is to investigate the behaviour of
super-Brownian motion among random obstacles, when the density of
these obstacles grows to infinity but their diameter shrinks to zero
in an appropriate manner. More precisely, let us fix $d\geq 2$ and a
domain $D$ of $\mathbb{R}^d$, and let $c:\mathbb{R}^d \rightarrow
[0,\infty)$ be a bounded measurable function. For every $\e \in
(0,\frac{1}{2})$, let us define an obstacle configuration by
$$\Gamma_{\e}=\bigcup_{x\in \mathcal{P}^{\e}}\overline{B}(x,\e),$$
where $\mathcal{P}^{\e}$ is a Poisson point process on
$\mathbb{R}^d$ with intensity $\log (\e^{-1})c(x)dx$ if $d=2$ and
$\e^{2-d}c(x)dx$ if $d\geq 3$, and $\overline{B}(x,\e)$ denotes the
closed ball of radius $\e$ centered at $x$. This Poisson point
process is defined on a probability space
$(\mathbf{\Omega},\mathbf{\mathcal{F}},\mathbf{P})$. On a different
probability space, let us also consider a superprocess $\{X^{\e}_t,\
t\in [0,\infty)\}$ with critical branching mechanism $\psi(u)=u^2$
and underlying spatial motion given by Brownian motion killed when
entering $D^c \cup \Gamma_{\e}$. Thus, for each $\e$, the
superprocess $X^{\e}$ can be seen as evolving in a random medium
given by $\Gamma_{\e}$. A realization of $\{\Gamma_{\e},\ \e\in
(0,1/2)\}$ will be called an \emph{environment}.

We wish to understand the behaviour of $X^{\e}$ when $\e$ tends to
zero. As in most works about random media, two points of view can be
adopted : either we fix an environment (quenched approach), or we
average over the possible realizations of
$\bigcup_{\e>0}\Gamma_{\e}$ (annealed approach). Although the
results of this paper are set in the quenched framework, the main
ingredients of their proofs are ``annealed-type'' calculations.
Moreover, the latter approach is also useful in obtaining a better
understanding of where the scaling comes from and of what the
limiting process might be. To simplify the analysis, let us first
assume that $D=\mathbb{R}^d$ and let us consider a single Brownian
motion $\xi$, independent of the obstacles. Denote by $\mathrm{P}_x$
the probability measure under which $\xi$ starts from $x$. Let us
define the random time $T_{\e}$ as the entrance time of $\xi$ into
the set $\Gamma_{\e}$, that is
$$ T_{\e}:=\inf \{t\geq 0: \xi_t \in \Gamma_{\e}\}.$$
In addition, for all $0\leq s\leq t$, we denote by $S_{\e}(s,t)$ the Wiener sausage of radius
$\e$ along the time interval $[s,t]$, defined as
$$S_{\e}(s,t)=\{y\in \mathbb{R}^d: \inf_{s\leq r\leq t}|\xi_r-y|\leq
\e\}= \bigcup_{r\in [s,t]}\big(\xi_r+\overline{B}(0,\e)\big).$$ The
probability that the Brownian motion $\xi$ hits $\Gamma_{\e}$ before
time $t$ is equal to the probability that the centre of one of the
obstacles lies in $S_{\e}(0,t)$. These centres are given by the
Poisson point process $\mathcal{P}^{\e}$ and so, by averaging over
the random obstacles and using Fubini's theorem, we obtain
\begin{equation}\label{eq intro}\mathbf{E}\big[\mathrm{P}_{0}[T_{\e} > t]\big]=
\mathrm{E}_0\big[\mathbf{P}[\mathcal{P}^{\e}
\cap S_{\e}(0,t)=\emptyset]\big] = \mathrm{E}_0\Big[\exp
-s_d(\e)\int_{S_{\e}(0,t)}c(x)dx\Big],
\end{equation}
where
$$s_d(\e)=\left\{\begin{array}{ll} \log (\e^{-1}) & \mathrm{if\ }d=2, \\
\e^{2-d} & \mathrm{if\ }d\geq 3. \end{array} \right. $$ In the case
$c=\nu$, the integral in (\ref{eq intro}) is just $\nu$ times the
volume $\lambda\big(S_{\e}(0,t)\big)$ of the Wiener sausage, whose
asymptotics have been well studied owing to their connections with
physical problems (see e.g. the introduction of \cite{SPI1964},
\cite{KL1974} or \cite{DV1975}). Note that the large-$t$ asymptotics
of $\lambda\big(S_{\e}(0,t)\big)$ are essentially equivalent to its
small-$\e$ asymptotics thanks to the equality in law :
$$\lambda\big(S_1(0,t)\big) \stackrel{(d)}{=} t^{d/2}\lambda\big(S_{t^{-1/2}}(0,1)\big).$$

A classical result of Kesten, Spitzer and Whitman (cf.
\cite{IM1965}, p.253) states that, if $d\geq 3$,
\begin{equation}\label{conv as intro}\lim_{\e \rightarrow 0}s_d(\e)\lambda\big(S_{\e}(0,t)\big)=k_d
t \qquad \mathrm{a.s.}, \end{equation} where
$k_d=(d-2)\pi^{d/2}/\Gamma(d/2)$ ($k_3=2\pi$) is the Newtonian
capacity of the unit ball. The Kesten-Spitzer-Whitman convergence
result was in fact stated for the large-time asymptotics of
$\lambda\big(S_{\e}(0,t)\big)$, but a scaling argument gives the
previous statement, at least in the sense of convergence in
probability. The convergence in (\ref{conv as intro}) also holds if
$d=2$ (see \cite{LEG1986}), with $k_2=\pi$.

It is not hard to deduce from the preceding result that, at least
when the function $c$ is continuous,
\begin{equation}\label{conv as intro c}\lim_{\e\rightarrow
0}s_d(\e)\int_{S_{\e}(0,t)}c(y)dy=k_d \int_0^tc(\xi_s)ds \qquad a.s.
\end{equation}
It then follows from (\ref{eq intro}) that
$$\lim_{\e\rightarrow 0}\mathbf{E}\big[\mathrm{P}_0[T_{\e}>
t]\big]= \mathrm{E}_0\Big[\exp -k_d\int_0^t c(\xi_s)ds\Big].$$ This
argument, which is due to Kac \cite{KAC1974}, can be interpreted in
the following way. When $\e$ tends to zero, the obstacles become
dense in $\mathbb{R}^d$ (at least if the function $c$ is everywhere
positive), and the Brownian motion $\xi_t$ gets absorbed in the
obstacles at rate $k_d c(\xi_t)$.

Back to our initial problem about killed superprocesses, the result
for a single Brownian particle suggests that the sequence $X^{\e}$
should converge to the superprocess $X^{*}$ with branching mechanism
$\psi(u,x)=u^2+k_d c(x) u$ and underlying spatial motion given by
Brownian motion. We shall establish in this work that the
distribution of $X^{*}$ is, indeed, the limit of the distribution of
$X^{\e}$ as $\e$ tends to $0$, in $\mathbf{P}$-probability. Here,
the distribution of $X^{\e}$ is a probability measure on the
Skorokhod space $D_{\mathcal{M}_f(\mathbb{R}^d)}([0,\infty))$ of all
c\`adl\`ag paths with values in $\mathcal{M}_f(\mathbb{R}^d)$ (the
space of all finite measures on $\mathbb{R}^d$) and the preceding
limit is in the sense of weak convergence. A stronger statement can
be made, but only for subsequences: if the sequence $\e_n$ decreases
to $0$ fast enough,
$$X^{\e_n} \stackrel{(d)}{\rightarrow} X^* \qquad \mathrm{as\ }n\rightarrow \infty,\ \mathbf{P}\mathrm{-a.s.}$$
Here, $\stackrel{(d)}{\rightarrow}$ denotes convergence in
distribution. Let us emphasize the meaning of this result : except
for a set of zero $\mathbf{P}$-measure, if we fix an environment,
then the sequence of superprocesses $X^{\e_n}$ evolving among these
fixed obstacles converges in law to $X^*$. Theorem \ref{main
theorem} and Corollary \ref{corollaire} are stated in a more general
setting, allowing the superprocesses to live only within a domain
$D$ of $\mathbb{R}^d$.

The question we address in this paper was motivated by analogous
works on Brownian motion. An extensive literature is already
available on this topic, reviewed for example in \cite{SZN1998}.
Owing to the well-known properties of Poisson point processes, they
seem to be a natural way to encode traps and have been frequently
exploited in investigations of the behaviour of Brownian motion
moving among ``hard'' obstacles, where the particle is killed
instantaneously when hitting an obstacle as described above, or
among ``soft'' obstacles, within which the Brownian particle is
killed at a certain rate. Our approach is close to ideas developed
by Kac in \cite{KAC1974}, whose probabilistic method differs from
the analytic method used by Papanicolaou and Varadhan \cite{PV1980}
in a similar context. Both derive the convergence in the
$L^2(\mathbf{P})$-norm of the semigroup of Brownian motion among
random obstacles when the number of obstacles tends to infinity but
their diameters tend to $0$ (recall that $\mathbf{P}$ denotes the
probability measure on the space where the obstacles are defined).
Subsequently, Brownian motion among traps was studied in different
settings, in particular by Sznitman, who devised the powerful method
of enlargement of obstacles (see \cite{SZN1998}).

The problem of super-Brownian motion or branching Brownian motion
among random obstacles was addressed recently by Engl\"ander in
\cite{ENG2000}, \cite{EH2002} and \cite{ENG2007}, the latter paper
dealing with soft obstacles. However, Engl\"ander considers the
supercritical case (instead of critical super-Brownian motion as we
do) and keeps the sizes of obstacles fixed. Within the obstacles, a
particle does not die but branches at a slower rate. His interest is
in the long-time asymptotics of the process and, in particular, the
survival probability and the growth rate of the support. His
techniques are mostly analytic, in contrast with the probabilistic
tools of the present work.

\subsection{Statement of the main result}
Let us first introduce some notation and construct the sequence of
superprocesses $X^{\e}$ from the historical superprocess
corresponding to a super-Brownian motion on $\mathbb{R}^d$,
independent of the obstacles. We refer to \cite{DP1991} for more
details on historical superprocesses and their applications. If $E$
is a topological space, $\mathcal{M}_f(E)$ stands for the space of
all finite Borel measures on $E$.

The (Brownian) historical superprocess can be defined as follows.
Let $\mathcal{W}$ be the set of all finite continuous paths in
$\mathbb{R}^d$, and note that $\mathbb{R}^d$ can be viewed as a
subset of $\mathcal{W}$ by identifying $x$ with the path of length
zero and initial point $x$. Then, let $\tilde{\xi}$ be the
continuous Markov process in $\mathcal{W}$ whose transition kernel
is described as follows: If $\tilde{\xi}_0=(w(r),\ 0\leq r\leq s)\in
\mathcal{W}$, the law of $\tilde{\xi}_t$ is the law under
$\mathrm{P}_{w(s)}$ of the concatenation of the paths $(w(r),\ 0\leq
r\leq s)$ and $(\xi_r,\ 0\leq r\leq t)$. The historical superprocess
$H$ is defined as the superprocess on $\mathcal{W}$ with branching
mechanism $\psi(u)=u^2$ and underlying spatial motion given by
$\tilde{\xi}$. Thus, $H$ takes values in
$\mathcal{M}_f(\mathcal{W})$. The super-Brownian motion $X^0$ on
$\mathbb{R}^d$, starting at $\mu \in\mathcal{M}_f(\mathbb{R}^d)$,
can then be recovered from the historical superprocess starting at
$\mu$ (which is viewed as a finite measure on the paths of length
zero) through the formula
$$\langle X^0_t,f \rangle =
\int_{\mathcal{W}}H_t(dw)f\big(w(t)\big)$$ for all $f$ bounded and
measurable and all $t\geq 0$. Here, $\langle \nu,f \rangle$ denotes
the integral of $f$ against the measure $\nu$.

We exploit this correspondence between the historical superprocess
and super-Brownian motion further to construct the sequence of
killed superprocesses which is of interest in this work. Let $E$ be
an open subset of $\mathbb{R}^d$, and recall the definition of the
obstacle configuration $\Gamma_{\e}$. For every $\e>0$, the
superprocess $\{X^{\e,E}_t, t\in [0,\infty)\}$ is defined from the
historical superprocess $H$ via the formula
$$\langle X^{\e,E}_t,f \rangle = \int_{\mathcal{W}}H_t(dw)f\big(w(t)\big)
\ind_{\{\forall s\in [0,t],\ w(s)\in E\cap \Gamma_{\e}^c\}},$$ for
all $f$ bounded and measurable, and all $t\geq 0$. It is
straightforward to verify that $X^{\e,E}$ is itself a super-Brownian
motion with critical branching mechanism $\psi(u)=u^2$ and
underlying spatial motion given by Brownian motion killed when
entering $E^c \cup \Gamma_{\e}$. Furthermore, $X^{\e,E}_0$ is the
restriction of $\mu$ to $E\cap \Gamma_{\e}^c$.

Recall that we defined $k_2=\pi$ and $k_d=
\frac{d-2}{\Gamma(d/2)}\pi^{d/2}$ for $d\geq 3$. We also introduce
another superprocess $X^{*,E}$, with branching mechanism
$\psi(u,x)=u^2+k_d c(x) u$ and underlying spatial motion given by
Brownian motion killed when it exits $E$.

In practice, $E$ will be either $D$ or a bounded open subset of $D$.
When there is no ambiguity, we shall suppress the dependence on $E$
in the notation. We choose a sequence $\e_n$ such that $\sum_n|\log
\e_n|^{-1}<\infty$ if $d=2$, and $\sum_n \e_n |\log\e_n |<\infty$ if
$d\geq 3$. For instance, we may fix $\alpha>1$ and set
$\e_n=\exp(-n^{\alpha})$ if $d=2$ and $\e_n=n^{-\alpha}$ if $d\geq
3$.

We will use the following notation.
\begin{itemize}
\item $\mathbb{P}_{\mu}$ is the (quenched) probability measure
under which $H$ starts at $\mu \in \mathcal{M}_f(\mathbb{R}^d)
\subset \mathcal{M}_f(\mathcal{W})$. By the preceding
correspondence, each superprocess $X^{\e,E}$ then starts under
$\mathbb{P}_{\mu}$ from the restriction of $\mu$ to $E\cap
\Gamma_{\e}^c$. It will be convenient to assume that $X^{*,E}$ is
also defined under $\mathbb{P}_{\mu}$ and starts from the
restriction of $\mu$ to $E$.
\item To simplify notation, $X^{(n),E}$ will be a shorthand for the killed superprocess with parameter
$\e_n$, and $\mathbb{P}^{(n),E}_{\mu}$ will be its law under
$\mathbb{P}_{\mu}$. Likewise, $\mathbb{P}^{\e,E}_{\mu}$ (resp.
$\mathbb{P}^{*,E}_{\mu}$) will be the law of $X^{\e,E}$ (resp.
$X^{*,E}$) under $\mathbb{P}_{\mu}$.
\item For all $t\geq 0$ and $x\in \mathbb{R}^d$, $\mathrm{P}_{t,x}$ will be a probability measure under
which a Brownian motion $\xi$ on $\mathbb{R}^d$, independent of the
obstacles, starts from $x$ at time $t$.
\item $T^E:=\inf \{t\geq 0: \xi_t \in E^c\}$, $T_{\e}:=\inf \{t\geq 0: \xi_t \in
\Gamma_{\e}\}$ and $T_{(n)}=T_{\e_n}$.
\end{itemize}

We can now state our main result.
\begin{theo}\label{main theorem} For every $\mu \in \mathcal{M}_f(D)$, $\mathbf{P}$-a.s.
$$\mathbb{P}^{(n),D}_{\mu} \Rightarrow \mathbb{P}^{*,D}_{\mu} \mathrm{\ as\ } n\rightarrow \infty,$$
where the symbol $\Rightarrow$ refers to the weak convergence of
probability measures.
\end{theo}

As an immediate corollary, we also have :
\begin{coro}\label{corollaire}
For every $\mu \in \mathcal{M}_f(D)$, the sequence
$\mathbb{P}^{\e,D}_{\mu}$ converges in $\mathbf{P}$-probability to
$\mathbb{P}^{*,D}_{\mu}$ as $\e$ tends to zero. In other words, for
every $\delta>0$, there exists $\e_0>0$ such that for all $0<\e\leq
\e_0$,
$$\mathbf{P}\left[d\big(\mathbb{P}^{\e,D}_{\mu},\mathbb{P}^{*,D}_{\mu}\big)> \delta
\right]< \delta,$$ where $d$ is the Prohorov metric on
$\mathcal{M}_1(D_{\mathcal{M}_f(D)}[0,\infty))$ (here,
$\mathcal{M}_1(D_{\mathcal{M}_f(D)}[0,\infty))$ is the space of all
probability measures on $D_{\mathcal{M}_f(D)}[0,\infty)$).
\end{coro}

The rest of the paper is devoted to the proofs of Theorem \ref{main
theorem} and Corollary \ref{corollaire}. In Section 2, we prove
certain estimates for the rate of convergence in (\ref{conv as intro
c}), which are of independent interest. These estimates are a key
ingredient of the proof of Lemma \ref{borne essentielle} in Section
3. Then, we fix a bounded open subset $B$ of $D$ and prove the
almost sure convergence of the distribution of $X^{(n),B}$ in two
steps. First, we show in Section $3$ that to each $k-$tuple
$(t_1,\ldots,t_k)$, there corresponds a set of $\mathbf{P}$-measure
zero outside which $(X_{t_1}^{(n),B},\ldots, X_{t_k}^{(n),B})_{n\geq
1}$ converges in law to $(X_{t_1}^{*,B}, \ldots, X_{t_k}^{*,B})$.
Second, we prove in Section $4$ that, with $\mathbf{P}$-probability
$1$, the sequence of superprocesses $X^{(n),B}$ is tight in
$D_{\mathcal{M}_f(D)}[0,\infty)$. In Section $5$, we complete the
proof for a general domain $D$. Starting with a bounded subset of
$D$ is required for technical reasons, to ensure the finiteness of
certain integrals which appear in the proof.

\section{Some estimates for the Wiener sausage}
Let us define the set $\mathcal{B}_1$ as the set of all bounded
Borel measurable functions $c$ on $\mathbb{R}^d$ such that
$\|c\|\leq 1$, where $\|c\|$ denotes the supremum norm of $c$. We
have the following result (we write $E_x$ for $E_{0,x}$ in the rest
of the section):
\begin{prop}\label{prop section 2}For every $t\geq 0$, there exists a constant $C=C(t)$ such
that for every $\e\in (0,\frac{1}{2}]$, if $d=2$,
$$\sup_{c \in \mathcal{B}_1}\sup_{x\in \mathbb{R}^2}\mathrm{E}_x\Big[\Big(|\log\e|\int_{S_{\e}(0,t)}c(y)dy -
\pi \int_0^t c(\xi_s)ds\Big)^2\Big]\leq \frac{C}{|\log\e|^2},$$ and
if $d\geq 3$,
$$\sup_{c \in \mathcal{B}_1}\sup_{x\in \mathbb{R}^d}\mathrm{E}_x\Big[\Big(\e^{2-d}\int_{S_{\e}(0,t)}c(y)dy -
k_d \int_0^t c(\xi_s)ds\Big)^2\Big]\leq C\e^2 |\log\e|^2.$$
\end{prop}

\begin{rema}In the case $c=1$, the bounds of Proposition \ref{prop section
2} follow from the known results for the fluctuations of the volume
of the Wiener sausage \cite{LEG1988}. However, it does not seem easy
to derive Proposition \ref{prop section 2} from the special case
$c=1$. Note that the latter case suggests that the bound $C\e^2|\log
\e|^2$ could be replaced by $C\e^2|\log \e|$ if $d=3$ and by $C\e^2$
if $d\geq 4$. These refinements will not be needed in our
applications.
\end{rema}

\noindent\textbf{Proof of Proposition \ref{prop section 2} for
$d\geq 3$:} To simplify notation, we prove the desired bound only
for $t=1$. A scaling argument then gives the result for any $t\geq
0$. Let us set
$$h(\e)=\sup_{c \in \mathcal{B}_1}\sup_{x\in \mathbb{R}^d}\mathrm{E}_x\Big[\Big(\e^{2-d}\int_{S_{\e}(0,1)}c(y)dy -
k_d \int_0^1 c(\xi_s)ds\Big)^2\Big].$$ As a first step, let us
notice that
$$\int_{S_{\e}(0,1)}c(y)dy = \int_{S_{\e}(0,1/2)}c(y)dy + \int_{S_{\e}(1/2,1)}c(y)dy -
\int_{S_{\e}(0,1/2)\cap S_{\e}(1/2,1)}c(y)dy.$$
Also,
$$\e^{2-d}\int_{S_{\e}(0,1/2)}c(y)dy - k_d \int_0^{1/2}c(\xi_s)ds = \e^{2-d}2^{-d/2}\int_{\tilde{S}_{\e\sqrt{2}}(0,1)}
c\Big(\frac{z}{\sqrt{2}}\Big)dz - \frac{k_d}{2}
\int_0^{1}c\Big(\frac{\tilde{\xi}_s}{\sqrt{2}}\Big)ds, $$ where
$\tilde{\xi}_s = \sqrt{2}\ \xi_{s/2}$ for all $s\geq 0$ and
$\tilde{S}_{\e}(a,b)$ is the Wiener sausage associated to
$\tilde{\xi}$. Since the function $\tilde{c}(z)=
c\big(\frac{z}{\sqrt{2}}\big)$ also belongs to $\mathcal{B}_1$, we
obtain that
$$\mathrm{E}_x\Big[\Big(\e^{2-d}\int_{S_{\e}(0,1/2)}c(y)dy - k_d \int_0^{1/2}c(\xi_s)ds\Big)^2\Big] \leq
\frac{1}{4}h(\e \sqrt{2}).$$ Likewise, using the Markov property at
time $\frac{1}{2}$ and the preceding argument, we have
$$\mathrm{E}_x\Big[\Big(\e^{2-d}\int_{S_{\e}(1/2,1)}c(y)dy - k_d \int_{1/2}^1 c(\xi_s)ds\Big)^2\Big] \leq
\frac{1}{4}h(\e \sqrt{2}).$$ On the other hand, we have
$\lambda\big(S_{\e}(0,1/2)\cap
S_{\e}(1/2,1)\big)=\lambda\big(S'_{\e}(0,1/2)\cap
S_{\e}''(0,1/2)\big)$, where $\xi_t'=\xi_{1/2-t}-\xi_{1/2}$ and
$\xi_t''=\xi_{1/2+t}-\xi_{1/2}$ for every $t\in [0,1/2]$, and
$S'_{\e}(0,1/2)$, resp. $S_{\e}''(0,1/2)$, denotes the Wiener
sausage with radius $\e$ associated to $\xi'$, resp. $\xi''$, along
the time interval $[0,1/2]$. Since $\xi'$ and $\xi''$ are
independent Brownian motions, we can use the following consequence
of Corollary 3-2 in \cite{LEG1986}, and of \cite{LEG1988}, p.1012:
There exists a constant $K_1(d)>0$ such that for every $\e \in
(0,1/2]$
$$\mathrm{E}\Big[\lambda\Big(S_{\e}(0,1/2)\cap S_{\e}(1/2,1)\Big)^2\Big]\leq
\left\{ \begin{array}{ll}K_1\e^4,& \quad d=3 \\ K_1 \e^8 |\log\e|^2, &\quad d=4 \\
K_1 \e^{2d}, &\quad d\geq 5. \end{array}\right.
$$
Coming back to the definition of $h(\e)$, and using the triangle
inequality in $L^2$, the fact that
$$\mathrm{E}_x\Big[\Big(\int_{S_{\e}(0,1/2)\cap
S_{\e}(1/2,1)}c(y)dy\Big)^2\Big]\leq
\mathrm{E}\Big[\lambda\Big(S_{\e}(0,1/2)\cap
S_{\e}(1/2,1)\Big)^2\Big]$$ and the preceding inequalities, we
obtain
\begin{eqnarray}h(\e) & \leq & \sup_{c \in \mathcal{B}_1}\sup_{x\in \mathbb{R}^d}
\bigg\{\mathrm{E}_x\Big[\Big(\e^{2-d}\int_{S_{\e}(0,1/2)}c(y)dy +
\e^{2-d}\int_{S_{\e}(1/2,1)}c(y)dy-k_d \int_0^1
c(\xi_s)ds\Big)^2\Big]^{1/2}\nonumber
\\ & & +\mathrm{E}_x\Big[\e^{4-2d}\Big(\int_{S_{\e}(0,1/2)\cap
S_{\e}(1/2,1)}c(y)dy\Big)^2\Big]^{1/2} \bigg\}^2 \nonumber\\
& \leq &  \Big\{\Big(\frac{1}{2}h(\e \sqrt{2})+
2u(\e)\Big)^{1/2}+K_1' \psi_d(\e)\Big\}^2, \label{eq1 section2}
\end{eqnarray} where $\psi_d(\e)=\e$ (resp. $\e^2|\log\e|$,
resp. $\e^2$) if $d=3$ (resp. $d=4$, resp. $d\geq 5$) and
$$u(\e)=\sup_{c\in \mathcal{B}_1}\sup_{x\in \mathbb{R}^d}
\Big|\mathrm{E}_x\Big[\Big(\e^{2-d}\int_{S_{\e}(0,1/2)}c(y)dy - k_d
\int_0^{1/2}c(\xi_s)ds\Big)\Big(\e^{2-d}\int_{S_{\e}(1/2,1)}c(y)dy -
k_d \int_{1/2}^1 c(\xi_s)ds\Big)\Big]\Big|.$$ Applying the Markov
property at time $\frac{1}{2}$, we have
$$u(\e)=\sup_{c\in \mathcal{B}_1}\sup_{x\in \mathbb{R}^d}
\Big|\mathrm{E}_x\Big[\Big(\e^{2-d}\int_{S_{\e}(0,1/2)}c(y)dy - k_d
\int_0^{1/2}c(\xi_s)ds\Big)v(\e,\xi_{1/2})\Big]\Big|,$$ where
$$v(\e,z)=\mathrm{E}_z\Big[\e^{2-d}\int_{S_{\e}(0,1/2)}c(y)dy - k_d
\int_0^{1/2}c(\xi_s)ds\Big].$$
We now use the following lemma.
\begin{lemm}\label{borne v,3}There exists a constant $K_2>0$ such
that for all $z\in \mathbb{R}^d$, $\e \in (0,\frac{1}{2}]$ and $c
\in \mathcal{B}_1$
$$|v(\e,z)|\leq K_2\ \e.$$
\end{lemm}
We postpone the proof of Lemma \ref{borne v,3} and complete the case
$d\geq 3$ of the Proposition. By Lemma \ref{borne v,3}, we have
$$
|u(\e)|\leq K_2\ \e \sup_{c\in \mathcal{B}_1}\sup_{x\in
\mathbb{R}^d}
\mathrm{E}_x\Big[\Big(\e^{2-d}\int_{S_{\e}(0,1/2)}c(y)dy - k_d
\int_0^{1/2}c(\xi_s)ds\Big)^2\Big]^{1/2}\leq \frac{K_2}{2}\ \e h(\e
\sqrt{2})^{1/2}.
$$
From (\ref{eq1 section2}), we obtain for every $\e \in
(0,\frac{1}{2}]$
$$h(\e)\leq \Big(\Big(\frac{1}{2}h(\e \sqrt{2})+ K_2\e h(\e \sqrt{2})^{1/2}\Big)^{1/2}+
K_1'\psi_d(\e) \Big)^2.$$ Let us set $g(\e)=\e^{-1}h(\e)^{1/2}$. We
thus have for $\e\in (0,\frac{1}{2}]$:
\begin{equation}\label{eqn g}g(\e)\leq
\big(g(\e\sqrt{2})^2+\sqrt{2}K_2g(\e\sqrt{2})\big)^{1/2}
+K_1'\e^{-1}\psi_d(\e).\end{equation} Fix $r\in (1/4,1/2]$ and set
$u_n=g(r2^{-n/2})$ for every integer $n\geq 0$. Rewriting (\ref{eqn
g}) in terms of $u_n$ and noting that $\e^{-1}\psi_d(\e)=1$ if $d=3$
and $\e^{-1}\psi_d(\e)=o(1)$ as $\e\rightarrow 0$ if $d\geq 4$, we
obtain for a constant $K_1''>0$ (independent of $n$)
$$u_{n+1}\leq (u_n^2+ \sqrt{2}K_2u_n)^{1/2}+ K_1'' =
u_n\Big(1+\frac{\sqrt{2}K_2}{u_n}\Big)^{1/2}+K_1''\leq u_n
+\frac{\sqrt{2}K_2}{2}+K_1''. $$ It follows that $u_n\leq u_0
+n\big(K_2 2^{-1/2}+K_1''\big)$ for every $n\geq 0$, from which we
can conclude that there exists a constant $K_3$ such that for all
$\e\in (0,1/2]$,
$$g(\e)\leq K_3 |\log \e|$$
and thus $$h(\e)\leq K_3^2\e^2|\log \e|^2.$$
\begin{flushright}
$\Box$ \end{flushright}

\noindent \textbf{Proof of Lemma \ref{borne v,3}:} We may assume
that $z=0$, and we fix $c\in \mathcal{B}_1$ (the constant $K_2$ will
not depend on $c$). First, we have
\begin{equation}\label{eq2 section2}\mathrm{E}_0\Big[k_d \int_0^{1/2}c(\xi_s)ds\Big]= k_d
\int_{\mathbb{R}^d}dy\ c(y) \int_0^{1/2} \frac{ds}{(2\pi
s)^{d/2}}\exp\Big(-\frac{|y|^2}{2s}\Big).\end{equation} Let us
define the random times $\tau_{\e}(y)$ and $L_{\e}(y)$ for all
$\e>0$ and $y\in \mathbb{R}^d$ by
\begin{eqnarray*}\tau_{\e}(y) &=& \inf\{t\geq 0:\ |\xi_t-y|\leq \e\}, \\
L_{\e}(y) &=& \sup\{t\geq 0:\ |\xi_t-y|\leq \e\},
\end{eqnarray*}
with the conventions that $\inf \emptyset = +\infty$ and $\sup
\emptyset = 0$. We thus have
\begin{eqnarray}\mathrm{E}_0\Big[\int_{S_{\e}(0,\frac{1}{2})}c(y)dy\Big]&=& \int_{\mathbb{R}^d}dy\ c(y)
\mathrm{P}_0\Big[\tau_{\e}(y)\leq \frac{1}{2}\Big]\label{eq3 section2}\\
&=& \int_{\mathbb{R}^d}dy\ c(y) \mathrm{P}_0\Big[0<L_{\e}(y)\leq
\frac{1}{2}\Big]+ \int_{\mathbb{R}^d}dy\ c(y)
\mathrm{P}_0\Big[\tau_{\e}(y)\leq
\frac{1}{2}<L_{\e}(y)\Big].\nonumber
\end{eqnarray}
On the one hand, $$\Big|\int dy\ c(y)
\mathrm{P}_0\Big[\tau_{\e}(y)\leq \frac{1}{2}<L_{\e}(y)\Big]\Big|
\leq \int dy\ \mathrm{P}_0\Big[\tau_{\e}(y)\leq \frac{1}{2}\leq
L_{\e}(y)\Big].$$ We have \begin{eqnarray*}\int dy\
\mathrm{P}_0\Big[\tau_{\e}(y)\leq \frac{1}{2}\leq L_{\e}(y)\Big]
&=& \mathrm{E}_0\Big[\lambda\Big(S_{\e}\big(0,1/2\big)\cap S_{\e}\big(1/2,\infty\big)\Big)\Big] \\
&=& \mathrm{E}_0\Big[\lambda\Big(S_{\e}\big(0,1/2\big)\cap
S'_{\e}\big(0,\infty\big)\Big)\Big],
\end{eqnarray*}
where $S'_{\e}$ denotes the Wiener sausage associated to a Brownian
motion $\xi'$ independent of $\xi$ and also started from $0$ under
$\mathrm{P}_0$. If $d=3$, it is easily checked that
\begin{equation}\label{eq proof lemma 1}\mathrm{E}_0\Big[\lambda\Big(S_{\e}\big(0,1/2\big)\cap
S'_{\e}\big(0,\infty\big)\Big)\Big]=O(\e^2)\end{equation} (use the
fact that $\mathrm{P}_0\big[y\in
S'_{\e}(0,\infty)\big]=\frac{\e}{|y|}\wedge 1$, together with the
bound (3.d) in \cite{LEG1988}). If $d\geq 4$,
\begin{eqnarray}
\mathrm{E}_0\Big[\lambda\Big(S_{\e}\big(0,1/2\big)\cap
S'_{\e}\big(0,\infty\big)\Big)\Big]&=&
\e^d\mathrm{E}_0\Big[\lambda\Big(S_1\big(0,\e^{-2}/2\big)\cap
S_1'\big(0,\infty\big)\Big)\Big] \nonumber \\
&=& \left\{\begin{array}{ll}O(\e^4|\log\e|), &\quad \mathrm{if}\
d=4,\\
O(\e^d),&\quad \mathrm{if}\ d\geq 5\end{array}\right. \label{eq4
section2} \end{eqnarray} by \cite{LEG1988}, p.1010.

Let $\nu_{\e,y}(dz)$ denote the equilibrium measure of the ball
$\overline{B}(y,\e)$, that is the unique finite measure on the
sphere $\partial B(y,\e)$ such that for every $x$ with $|x-y|>\e$,
$$\mathrm{P}_x\big[\tau_{\e}(y)<\infty\big]=\int
\nu_{\e,y}(dz)G(z-x),$$ where $G(z)=\int_0^{\infty}(2\pi
s)^{-d/2}\exp(-|z|^2/2s)ds = c_d |y|^{2-d}$ is the Green function of
$d$-dimensional Brownian motion ($c_d$ is a constant depending only
on $d$). By a classical formula of probabilistic potential theory
(see \cite{PS1978}, p.61-62) we have
$$\mathrm{P}_0\Big[0<L_{\e}(y)\leq \frac{1}{2}\Big]= \int_0^{1/2} ds \int \nu_{\e,y}(dz)\frac{1}{(2\pi s)^{d/2}}
\exp\Big(-\frac{|z|^2}{2s}\Big).$$ It is well known that
$\nu_{\e,y}= k_d \e^{d-2} \pi_{\e,y}$, where $\pi_{\e,y}$ denotes
the uniform distribution on the sphere of radius $\e$ centered at
$y$. Recalling (\ref{eq2 section2}), (\ref{eq3 section2}), (\ref{eq
proof lemma 1}) and (\ref{eq4 section2}), we can write
\begin{eqnarray*}|v(\e,z)|&=& \bigg| k_d \int_{\mathbb{R}^d}dy\ c(y) \int_0^{1/2}
\frac{ds}{(2\pi s)^{d/2}}\bigg\{  \int \pi_{\e,y}(dz)
\big(e^{-|z|^2/2s}-e^{-|y|^2/2s} \big) \bigg\}
\bigg| +O(\phi_d(\e))\\
&\leq & k_d \int_{\mathbb{R}^d}dy \int_0^{1/2}\frac{ds}{(2\pi
s)^{d/2}}\bigg\{\int \pi_{\e,y}(dz) \Big|e^{-|z|^2/2s}-e^{-|y|^2/2s}
\Big| \bigg\}+ O(\phi_d(\e))
\end{eqnarray*}
where $\phi_d(\e)=\e$ (resp. $\e^2|\log\e|$, resp. $\e^2$) if $d=3$
(resp. $d=4$, resp. $d\geq 5$). It follows that
\setlength\arraycolsep{1pt}
\begin{eqnarray*}\int_{|y|\leq 10\e}&dy&\int_0^{1/2}\frac{ds}{(2\pi
s)^{d/2}}\int
\pi_{\e,y}(dz)\Big|\ e^{-|z|^2/2s}-e^{-|y|^2/2s}\Big|\\
&\leq & 2\int_{|z|\leq 11\e}dz\int_0^{1/2} \frac{ds}{(2\pi
s)^{d/2}}\ e^{-|z|^2/2s} \leq 2\int_{|z|\leq 11\e}dz\ G(z) =O(\e^2).
\end{eqnarray*}
On the other hand, we can find constants $C$ and $C'$ such that if
$|y|>10\e$ and $|z-y|=\e$,
$$\Big|\exp\Big(-\frac{|z|^2}{2s}\Big)-\exp\Big(-\frac{|y|^2}{2s}\Big)\Big|\leq
C\bigg|\frac{|z|^2-|y|^2}{s}\bigg|\exp\Big(-\frac{|y|^2}{4s}\Big)\leq
C'\e\frac{|y|}{s}\exp\Big(-\frac{|y|^2}{4s}\Big).$$ Thus, with a
constant $K$ which may vary from line to line, we have
\setlength\arraycolsep{1pt}
\begin{eqnarray*}
\int_{|y|>10\e}&dy& \int_0^{1/2}\frac{ds}{(2\pi s)^{d/2}}\int
\pi_{\e,y}(dz)\Big|\exp\Big(-\frac{|z|^2}{2s}\Big)-\exp\Big(-\frac{|y|^2}{2s}\Big)\Big|\\
&\leq& K\e \int_{|y|>10\e}dy\ |y|\int_0^{1/2}ds\
s^{-d/2-1}\exp\Big(-\frac{|y|^2}{4s}\Big) \\
&=& K\e \int_{|y|>10\e}dy\ |y|^{1-d}\int_0^{1/(2|y|^2)}ds'\
s'^{-d/2-1}e^{-1/4s'}\\
&\leq &K\e.
\end{eqnarray*}
Combining the above, the proof of Lemma \ref{borne v,3} is complete.

\bigskip \noindent\textbf{Proof of Proposition \ref{prop section 2}
for $d=2$:} Let us define
$$h(\e)=\sup_{c\in \mathcal{B}_1}\sup_{x\in \mathbb{R}^d} \mathrm{E}\Big[\Big(|\log \e|\int_{S_{\e}(0,1)}c(y)dy
-\pi \int_0^1 c(\xi_s)ds\Big)^2\Big].$$ By Corollary 3-2 in
\cite{LEG1986}, we have
$$\mathrm{E}\Big[\lambda\Big(S_{\e}\big(0,\frac{1}{2}\big)\cap S_{\e}\big(\frac{1}{2},1)\Big)^2\Big]
\leq \frac{K_1}{|\log \e|^4}.$$ The same technique as in the
previous case yields
\begin{equation} \label{borne h,2}h(\e)\leq \Big(\Big(\frac{1}{2}h(\e\sqrt{2})+h(\e\sqrt{2})^{1/2} \sup_{c
\in \mathcal{B}_1}\sup_{x\in \mathbb{R}^2}|v(\e,z)|\Big)^{1/2}+
\frac{\sqrt{K_1}}{|\log \e|} \Big)^2,\end{equation} where $$v(\e,z)=
\mathrm{E}_z\Big[|\log \e|\int_{S_{\e}(0,1/2)}c(y) dy - \pi
\int_0^{1/2}c(\xi_s)ds \Big].$$ We now use the following result,
whose rather technical proof is deferred to the Appendix : There
exists a constant $K_2$ such that, for $\e\in (0,1/2]$,
\begin{equation}\label{borne v,2}\sup_{c\in \mathcal{B}_1}\sup_{z\in \mathbb{R}^d}|v(\e,z)| \leq \frac{K_2}{|\log \e|}.
\end{equation}
Hence, if $g(\e)=|\log\e|h(\e)^{1/2}$, we have for $\e\in (0,1/2]$,
\begin{equation}\label{ineg g}g(\e)\leq \bigg(\frac{1}{2}\frac{(\log\e)^2}{(\log\e\sqrt{2})^2}\ g(\e\sqrt{2})^2
+K_2\frac{|\log\e|}{|\log\e\sqrt{2}|}\
g(\e\sqrt{2})\bigg)^{1/2}+\sqrt{K_1}.\end{equation} From (\ref{ineg
g}), we can use arguments similar to the case $d\geq 3$ to infer
that the function $g(\e)$ is bounded over $(0,1/2]$. Thus, there
exists a constant $K_3$ such that for all $\e \in (0,1/2]$,
$$h(\e)\leq \frac{K_3}{|\log \e|^2}.$$
\begin{flushright}$\Box$ \end{flushright}

\section{Almost sure convergence of the finite-dimensional distributions of $X^{(n),B}$}
In the following, we fix a bounded open subset $B$ of $D$ and
consider only the superprocesses killed outside $B$. We therefore
suppress the dependence on $B$ in the notation. In particular,
$T=T^{B}$.

The following proposition is the first step in the proof of Theorem
\ref{main theorem}.
\begin{prop}\textbf{(Convergence of the finite-dimensional distributions)}
Let $\mu \in \mathcal{M}_f(B)$, $p\in \mathbb{N}$ and $t_1< \ldots
<t_p \in [0,\infty)$. Then, under $\mathbb{P}_{\mu}$
$$(X_{t_1}^{(n)},\ldots, X_{t_p}^{(n)}) \stackrel{(d)}{\rightarrow} (X_{t_1}^*,\ldots,
X_{t_p}^*)$$ as $n\rightarrow \infty$, on a set of
$\mathbf{P}$-probability $1$. \label{conv_fini_dim}
\end{prop}

\paragraph{Proof of Proposition \ref{conv_fini_dim} :}
We fix an environment. Let $p \in \mathbb{N}$, $0\leq
t_1<\ldots<t_p$ and $f_1,\ldots,f_p \in
\mathcal{B}_{b+}(\mathbb{R}^d)$ be measurable, nonnegative and
bounded functions. In the following, we shall denote
$(t_1,\ldots,t_p)$ by $\mathbf{t}$ and $(f_1,\ldots,f_p)$ by
$\mathbf{f}$.

Let $\mu \in \mathcal{M}_f(B)$. Following the notation in
\cite{DYN1991}, we have :
\begin{eqnarray*}
\mathbb{E}_{\mu}\Big[\exp -\sum_{i=1}^p\langle X^{\e}_{t_i}, f_i\rangle\Big] &=& \exp -\langle \mu, w_0^{\e}\rangle,\\
\mathbb{E}_{\mu}\Big[\exp -\sum_{i=1}^p\langle X^{*}_{t_i},
f_i\rangle\Big] &=& \exp -\langle \mu, w_0^{*}\rangle,
\end{eqnarray*}
where $w^{\e}=(w^{\e}_t(x);\ t\geq 0,\ x\in B)$ and $w^*=(w^*_t(x);\
t\geq 0,\ x\in B)$ are the unique nonnegative solutions to the
following integral equations: for all $x\in B$ and $t\geq 0$,
\begin{eqnarray}w_t^{\e}(x)+ \mathrm{E}_{t,x}\Big[\int_t^{\infty}ds\ w_s^{\e}(\xi_s)^2
\ind_{\{s<T\wedge T_{\e}\}} \Big] & =&
\sum_{i=1}^p \mathrm{E}_{t,x}\Big[f_i(\xi_{t_i})\ind_{\{t_i<T\wedge T_{\e}\}}\Big], \label{caract_we} \\
w_t^{*}(x)+ \mathrm{E}_{t,x}\Big[\int_t^{\infty}ds\
\big(w_s^{*}(\xi_s)^2+k_dc(\xi_s)w_s^*(\xi_s)\big)\ \ind_{\{s<T\}}
\Big] & = & \sum_{i=1}^p \mathrm{E}_{t,x}\Big[f_i(\xi_{t_i})\
\ind_{\{t_i<T\}}\Big] \label{caract_w*},
\end{eqnarray}
where by convention $\mathrm{E}_{t,x}[f(\xi_s)]=0$ if $s<t$. By the
standard argument of the proof of the Feynman-Kac formula, the
integral equation (\ref{caract_w*}) for $w^*$ is equivalent to
\begin{equation}\label{caract_w* bis}w_t^{*}(x)+ \mathrm{E}_{t,x}\Big[\int_t^{\infty}ds\
w_s^{*}(\xi_s)^2 e^{-k_d\int_t^s c(\xi_u)du}\ \ind_{\{s<T\}} \Big] =
\sum_{i=1}^p
\mathrm{E}_{t,x}\Big[f_i(\xi_{t_i})e^{-k_d\int_t^{t_i}c(\xi_u)du}\
\ind_{\{t_i<T\}}\Big]. \end{equation}
The equivalence of the two
integral equations (\ref{caract_w*}) and (\ref{caract_w* bis})
corresponds to the well-known fact that super-Brownian motion with
branching mechanism $\psi(u,x)=u^2 +k_dc(x)u$ can also be
constructed as the superprocess with branching mechanism
$\psi(u)=u^2$ and underlying spatial motion given by Brownian motion
killed at rate $k_dc(x)$.

\begin{rema} Since $w^{\e}$ and
$w^{*}$ are nonnegative, (\ref{caract_we}) and (\ref{caract_w*})
imply that $w^{*}_t$ and $w^{\e}_t$ are equal to zero whenever $t>
t_p$ (recall that by convention, the right-hand side of
(\ref{caract_we}) or (\ref{caract_w*}) is zero when $t>t_p$).
Likewise, $w_t^{\e}(x)=0$ if $x\in \Gamma_{\e}\cap B$, for every
$t\geq 0$.
\end{rema}

By integrating over $B$ the difference between (\ref{caract_we}) and
(\ref{caract_w* bis}), we obtain : \setlength\arraycolsep{1pt}
\begin{eqnarray}
\int_{B}dx\ |w_t^{\e}(x)-w_t^{*}(x)| &\leq& \int_{B}dx\
\Big|\sum_{i=1}^p \mathrm{E}_{t,x}\Big[
f_i(\xi_{t_i})\ind_{\{t_i<T\}}
(\ind_{\{t_i<T_{\e}\}}-e^{-k_d\int_t^{t_i}c(\xi_u)du})\Big]\Big| \label{terme1}\\
& & +\int_{B}dx\ \Big|\mathrm{E}_{t,x}\Big[\int_t^{\infty}ds\
\ind_{\{s<T\}}(\ind_{\{s<T_{\e}\}}-e^{-k_d\int_t^s
c(\xi_u)du})w_s^{\e}(\xi_s)^2 \Big]\Big|
\nonumber\\
& & + \int_{B}dx\ \Big|\mathrm{E}_{t,x}\Big[\int_t^{\infty}ds\
e^{-k_d\int_t^s c(\xi_u)du}\ind_{\{s<T\}} (w_s^{*}(\xi_s)^2 -
w_s^{\e}(\xi_s)^2)\Big]\Big|.
 \nonumber
\end{eqnarray}
Let us start with the third term in the right-hand side of
(\ref{terme1}). The functions $w^{\e}$ and $w^*$ are bounded by
$C_{\mathbf{f}}:= \sum_{i=1}^p \|f_i\|$, hence bounding
$\ind_{\{s<T\}}$ by $\ind_{B}(\xi_s)$ and $e^{-k_d\int_t^s
c(\xi_u)du}$ by $1$ yields \setlength\arraycolsep{1pt}
\begin{eqnarray}
\int_{B}dx& & \Big|\mathrm{E}_{t,x}\Big[\int_t^{\infty}ds\
e^{-k_d\int_t^s c(\xi_u)du}\ind_{\{s<T\}} (w_s^{*}(\xi_s)^2 -
w_s^{\e}(\xi_s)^2)
\Big]\Big| \nonumber \\
& &\leq 2 C_{\mathbf{f}} \int_{B}dx\
\mathrm{E}_{t,x}\Big[\int_t^{\infty}ds\ \ind_{B}(\xi_s)
|w_s^{*}(\xi_s)-w_s^{\e}(\xi_s)|\Big]
\nonumber \\
& &=\ 2 C_{\mathbf{f}} \int_t^{\infty}ds\ \int_{B \times B}dx\ dz\ |w_s^{*}(z)-w_s^{\e}(z)|p_{s-t}(x,z) \nonumber\\
& &\leq 2 C_{\mathbf{f}} \int_t^{\infty}ds\ \int_{B}dz\
|w_s^{*}(z)-w_s^{\e}(z)|. \label{gronwall}
\end{eqnarray}
In the preceding estimates, $p_r(\cdot,\cdot)$ denotes the
transition density at time $r$ of $d$-dimensional Brownian motion.
The last inequality stems from the observation that
$\int_{B}p_{s-t}(x,z)dx = \int_{B}p_{s-t}(z,x)dx \leq 1$.

We next show that the first two terms of (\ref{terme1}) converge
towards $0$ $\mathbf{P}$-a.s. The key ingredient is the following
result:
\begin{lemm}\label{borne essentielle}
Let $t_1\in [0,\infty)$ and let $f \in
\mathcal{B}_{b+}(\mathbb{R}^d)$ be a bounded nonnegative measurable
function. Then, there exists a constant $K=K(c,t_1,d)$ such that,
for every $t\in [0,\infty)$, $x\in B$ and $\e\in (0,1/2)$, if $\
d=2$
$$
\mathbf{E}\Big[\mathrm{E}_{t,x}\Big[f(\xi_{t_1})\ind_{\{t_1<T\}}\big(\ind_{\{t_1<T_{\e}\}}-e^{-\pi\int_t^{t_1}
c(\xi_u)du}\big)\Big]^2\Big] \leq K \ \|f\|^2\ \frac{1}{|\log \
\e|},
$$
and if $\ d\geq 3$,
$$
\mathbf{E}\Big[\mathrm{E}_{t,x}\Big[f(\xi_{t_1})\ind_{\{t_1<T\}}\big(\ind_{\{t_1<T_{\e}\}}-e^{-k_d\int_t^{t_1}
c(\xi_u)du}\big)\Big]^2 \Big] \leq K\ \|f\|^2\ \e |\log\ \e|.$$
\end{lemm}

The proof of Lemma \ref{borne essentielle} is postponed until the
end of the section. Let us temporarily fix $t\in [0,t_p]$. Applying
the lemma with $\e=\e_n$, we obtain for every $\delta>0$ and every
$i\in \{1,\ldots,p\}$ \setlength\arraycolsep{1pt}
\begin{eqnarray*}\mathbf{P}\Big[\int_{B}&\Big|&\mathrm{E}_{t,x}\Big[f_i(\xi_{t_i})\ind_{\{t_i<T\}}
(\ind_{\{t_i<T_{(n)}\}} -\
e^{-k_d\int_t^{t_i}c(\xi_u)du})\Big]\Big|\ dx
  >\delta \Big] \\&& \leq \frac{1}{\delta^2}\ \mathbf{E}\Big[\Big(\int_{B}\Big|\mathrm{E}_{t,x}
  \Big[f_i(\xi_{t_i})\ind_{\{t_i<T\}}
  (\ind_{\{t_i<T_{(n)}\}}-e^{-k_d\int_t^{t_i}c(\xi_u)du})\Big]\Big|\ dx \Big)^2\Big] \\
  &&\leq \frac{\lambda(B)}{\delta^2}\int_{B}\mathbf{E}\Big[\mathrm{E}_{t,x}\Big[f_i(\xi_{t_i})\ind_{\{t_i<T\}}
  (\ind_{\{t_i<T_{(n)}\}}-e^{-k_d\int_t^{t_i}c(\xi_u)du})\Big]^2\Big]dx \\
  &&\leq \left\{\begin{array}{ll}
  \lambda(B)^2 K\
  \|f_i\|^2 \delta^{-2}\ |\log\e_n|^{-1} & \mathrm{if\ }d=2, \\
  \lambda(B)^2 K\ \|f_i\|^2 \delta^{-2}\ \e_n|\log\e_n| & \mathrm{\ if\ }d\geq 3,
  \end{array}\right.
\end{eqnarray*}
which is summable by our assumptions on $(\e_n)_{n\geq 1}$. Hence,
by the Borel-Cantelli lemma,
$$\mathbf{P}\mathrm{-a.s.,\
}\int_{B}\Big|\mathrm{E}_{t,x}\Big[f_i(\xi_{t_i})\ind_{\{t_i<T\}}(\ind_{\{t_i<T_{(n)}\}}-e^{-k_d\int_t^{t_i}
c(\xi_u)du})\Big]\Big|\ dx \rightarrow 0$$ as $n$ tends to infinity.
The first term of (\ref{terme1}) is bounded above by a finite sum of
such terms, therefore it converges to $0$ $\mathbf{P}$-a.s, for each
fixed $t\in [0,t_p]$.

Let us set
$$A_{\mathbf{f},\mathbf{t}}:=\left\{(\omega,t)\in \Omega\times
[0,t_p]:\int_{B}dx\ \Big|\mathrm{E}_{t,x}\Big[\sum_{i=1}^p
f_i(\xi_{t_i})\ind_{\{t_i<T\}}(\ind_{\{t_i<T_{(n)}\}}-e^{-k_d\int_t^{t_i}c(\xi_u)du})\Big]\Big|
\rightarrow 0 \right\}.$$ If $\lambda_1$ denotes the Lebesgue
measure on $\mathbb{R}$, we have by Fubini's theorem
$\mathbf{P}\otimes
\lambda_1(A_{\mathbf{f},\mathbf{t}}^c)=\int_0^{t_p} dt\
\mathbf{P}\big(\{\omega:(\omega,t)\in
A_{\mathbf{f},\mathbf{t}}^c\}\big)=0$, which gives $(i)$ in the
following lemma :
\begin{lemm}\label{exist_set}

(i) There exists a measurable subset
$\tilde{\Omega}_{\mathbf{f},\mathbf{t}}$ of $\mathbf{\Omega}$, with
$\mathbf{P}(\tilde{\Omega}_{\mathbf{f},\mathbf{t}})=0$, such that
for every $\omega~\in~\mathbf{\Omega} \backslash
\tilde{\Omega}_{\mathbf{f},\mathbf{t}}$,
$$\int_{B}\Big|\mathrm{E}_{t,x}\Big[\sum_{i=1}^p
f_i(\xi_{t_i})\ind_{\{t_i<T\}}(\ind_{\{t_i<T_{(n)}\}}-e^{-k_d\int_t^{t_i}c(\xi_u)du})\Big]\Big|\
dx \rightarrow 0 \quad \mathrm{as}\: n\rightarrow \infty $$ for all
$t\geq 0$, except for $t$ belonging to a Lebesgue null subset
$\tilde{T}_{\mathbf{f},\mathbf{t},\omega}$ of $\mathbb{R}_+$.

(ii) There exists also a measurable subset
$\hat{\Omega}_{\mathbf{f},\mathbf{t}}$ of $\mathbf{\Omega}$, with
$\mathbf{P}(\hat{\Omega}_{\mathbf{f},\mathbf{t}})=0$, such that for
every $\omega \in \mathbf{\Omega}\backslash
\hat{\Omega}_{\mathbf{f},\mathbf{t}}$,
$$\int_0^{\infty}ds\
\int_{B}dx\
\Big|\mathrm{E}_{0,x}\Big[\ind_{\{s<T\}}(\ind_{\{s<T_{(n)}\}}-e^{-k_d\int_0^sc(\xi_u)du})w^{(n)}_{s+t}(\xi_s)^2
\Big]\Big| \rightarrow 0 \qquad \mathrm{as}\: n\rightarrow \infty$$
for all $t\geq 0$, except on a Lebesgue null subset
$\hat{T}_{\mathbf{f},\mathbf{t},\omega}$ of $\mathbb{R}_+$. Here,
$w^{(n)}=w^{\e_n}$ is the function given by (\ref{caract_we})
corresponding to the superprocess $X^{(n)}$.

(iii) Finally, for all $x \in B$ there exists a negligible set
$\Omega_{\mathbf{f},\mathbf{t},0}(x)$ outside which
$$\Big|\mathrm{E}_{0,x}\Big[\sum_{i=1}^p f_i(\xi_{t_i})\ind_{\{t_i<T\}}
(\ind_{\{t_i<T_{(n)}\}}-e^{-k_d\int_0^{t_i}c(\xi_u)du})\Big]\Big| $$
and
$$ \Big|\mathrm{E}_{0,x}\Big[\int_0^{\infty}ds\
\ind_{\{s<T\}}(\ind_{\{s<T_{(n)}\}}-e^{-k_d\int_0^sc(\xi_u)du})w_s^{(n)}(\xi_s)^2
\Big]\Big|
$$ converge to $0$ as $n\rightarrow \infty$.
\end{lemm}
Both $(ii)$ and $(iii)$ can be obtained from Lemma \ref{borne
essentielle} in a way similar to the derivation of $(i)$. Note that
in $(ii)$, we may replace the integral over $[0,\infty)$ by the
integral over $[0,t_p]$ (since $w_r^{(n)}\equiv 0$ if $r\geq t_p$)
and that the functions $w_r^{(n)}$ are uniformly bounded by
$C_{\mathbf{f}}$.

The first term of the right-hand side of (\ref{terme1}), with
$\e=\e_n$, converges to $0$ as $n\rightarrow \infty$ provided that
$\omega \notin \tilde{\Omega}_{\mathbf{f},\mathbf{t}}$ and $t\notin
\tilde{T}_{\mathbf{f},\mathbf{t},\omega}$, by Lemma \ref{exist_set}
$(i)$. For the second term, we have \setlength\arraycolsep{1pt}
\begin{eqnarray}
\int_{B}dx\ \Big|\ \mathrm{E}_{t,x}&&\Big[\int_t^{\infty}ds\
\ind_{\{s<T\}}\big(\ind_{\{s<T_{(n)}\}}-e^{-k_d\int_t^sc(\xi_u)du}\big)w_s^{(n)}(\xi_s)^2
\Big]\Big| \label{borne_terme_2}\\
& & = \int_{B}dx\ \Big|\ \mathrm{E}_{0,x}\Big[\int_0^{\infty}ds\
\ind_{\{s<T\}}\big(\ind_{\{s<T_{(n)}\}}-e^{-k_d\int_0^sc(\xi_u)du}\big)w_{s+t}^{(n)}(\xi_s)^2
\Big]\Big| \nonumber \\
& & \leq \int_0^{\infty}ds\ \int_{B}dx\ \Big|\
\mathrm{E}_{0,x}\Big[\ind_{\{s<T\}}\big(\ind_{\{s<T_{(n)}\}}-e^{-k_d\int_0^sc(\xi_u)du}\big)w_{s+t}^{(n)}(\xi_s)^2
\Big]\Big|\ , \nonumber
\end{eqnarray}
which converges to $0$ as $n\rightarrow \infty$ by Lemma
\ref{exist_set} $(ii)$, if $\omega \notin
\hat{\Omega}_{\mathbf{f},\mathbf{t}}$ and $t \notin
\hat{T}_{\mathbf{f},\mathbf{t},\omega}$.

Finally, for $\omega \in \big(\tilde{\Omega}_{\mathbf{f},\mathbf{t}}
\cup \ \hat{\Omega}_{\mathbf{f},\mathbf{t}}\big)^c$ and $t\in
\big(\tilde{T}_{\mathbf{f},\mathbf{t},\omega} \cup\
\hat{T}_{\mathbf{f},\mathbf{t},\omega}\big)^c$, the first two terms
of the right-hand side of (\ref{terme1}) converge to $0$ as
$n\rightarrow \infty$. Recalling (\ref{gronwall}), we obtain
$$\int_{B}dx\ |w_t^{(n)}(x)-w_t^*(x)| \leq b_n(t)+2C_{\mathbf{f}} \int_t^{t_p}ds\ \int_{B}dz\ |w_t^{(n)}(z)-w_t^*(z)|,$$
where $b_n(t)\rightarrow 0$ as $n\rightarrow \infty$ provided
$\omega$ and $t$ are as above. Besides, for every $t$,
\begin{equation}\label{borne b} |b_n(t)|\leq
2\lambda(B)(1+t_p)(C_{\mathbf{f}}+C_{\mathbf{f}}^2).\end{equation}
Set for every $t\in [0,t_p]$,
$$G_n(t):= \int_{B}dx\ |w_{t_p-t}^{(n)}(x)-w_{t_p-t}^*(x)|.$$
Then, $$G_n(t) \leq b_n(t_p-t)+K_{\mathbf{f}} \int_0^t ds\ G_n(s),$$
where $K_{\mathbf{f}}:= 2C_{\mathbf{f}}$. By iterating this
inequality as in the proof of Gronwall's lemma, we obtain for all
$k\geq 1$, $n\geq 1$ and $t\in [0,t_p]$
\begin{eqnarray*}
G_n(t) &\leq& b_n(t_p-t)+
\sum_{i=0}^{k-2}K_{\mathbf{f}}^{i+1}\int_0^tds_1\int_0^{s_1}ds_2\ldots\int_0^{s_i}ds_{i+1}\
b_n(t_p-s_{i+1}) \\ & & \qquad \qquad \qquad + K_{\mathbf{f}}^k\int_0^tds_1\int_0^{s_1}ds_2\ldots\int_0^{s_{k-1}}ds_k G_n(s_k) \\
&\leq& b_n(t_p-t)+
\sum_{i=0}^{k-2}K_{\mathbf{f}}^{i+1}\int_0^tds_1\int_0^{s_1}ds_2\ldots\int_0^{s_i}ds_{i+1}\
b_n(t_p-s_{i+1}) + \lambda(B)K_{\mathbf{f}}^{k+1}\frac{t_p^k}{k!}\ .
\end{eqnarray*}
Fix $\e>0$ and let $k\geq 2$ be such that
$\lambda(B)K_{\mathbf{f}}^{k+1}\frac{t_p^k}{k!}\leq \frac{\e}{2}$.
For $\omega \notin \big(\tilde{\Omega}_{\mathbf{f},\mathbf{t}} \cup
\ \hat{\Omega}_{\mathbf{f},\mathbf{t}}\big)$, $b_n(r)$ converges to
$0$ as $n\rightarrow \infty$ except on a Lebesgue null set of values
of $r$, and thus by dominated convergence
$$\int_0^tds_1\int_0^{s_1}ds_2\ldots\int_0^{s_i}ds_{i+1}\
b_n(t_p-s_{i+1})\rightarrow 0$$ for every $t\in [0,t_p]$ and $i\in
\{0,\ldots,k-2\}$. In particular, for such $\omega$ and for every
$t\in [0,t_p]$, we have
$$K_{\mathbf{f}}^{i+1}\int_0^tds_1\int_0^{s_1}ds_2\ldots\int_0^{s_i}ds_{i+1}\
b_n(t_p-s_{i+1})\leq \frac{\e}{2k}$$ for all $n$ sufficiently large.
If moreover $t$ is such that $t_p-t\in
\Big(\tilde{T}_{\mathbf{f},\mathbf{t},\omega} \cup
\hat{T}_{\mathbf{f},\mathbf{t},\omega}\Big)^c$, then we have also
$b_n(t_p-t)\leq \frac{\e}{2k}$ if $n$ is large. Hence, we have
$G_n(t)\leq \e$ when $n$ is large. Since $\e$ was arbitrary, we can
conclude that for all $\omega$ and $t$ as specified above, $G_n(t)$
converges to $0$. Equivalently: for all $\omega \in
\big(\tilde{\Omega}_{\mathbf{f},\mathbf{t}} \cup
\hat{\Omega}_{\mathbf{f},\mathbf{t}}\big)^c$ and $t\in
\big(\tilde{T}_{\mathbf{f},\mathbf{t},\omega} \cup
\hat{T}_{\mathbf{f},\mathbf{t},\omega}\big)^c$,
\begin{equation}\label{convergence1}\lim_{n\rightarrow \infty}\int_{B}dx\
|w_t^{(n)}(x)-w_t^*(x)|= 0.
\end{equation}

We next consider the asymptotic behaviour of
$|w_0^{(n)}(x)-w_0^*(x)|$. In the same way as in (\ref{terme1}) but
now without integrating over $B$, we have for every $x\in B$
\begin{eqnarray}
|w_0^{(n)}(x)-w_0^{*}(x)| &\leq&
\Big|\mathrm{E}_{0,x}\Big[\sum_{i=1}^p
f_i(\xi_{t_i})\ind_{\{t_i<T\}}
(\ind_{\{t_i<T_{(n)}\}}-e^{-k_d\int_0^{t_i}c(\xi_u)du})\Big]\Big| \label{terme1bis}\\
& & +\ \Big|\mathrm{E}_{0,x}\Big[\int_0^{\infty}ds\
\ind_{\{s<T\}}(\ind_{\{s<T_{(n)}\}}-e^{-k_d\int_0^sc(\xi_u)du})w_s^{(n)}(\xi_s)^2 \Big]\Big| \nonumber \\
& & +\ 2C_{\mathbf{f}} \mathrm{E}_{0,x}\Big[\int_0^{\infty}ds\
\ind_{B}(\xi_s) |w_s^{*}(\xi_s) - w_s^{(n)}(\xi_s)|\Big]. \nonumber
\end{eqnarray}
Let us fix $x\in B$. By Lemma \ref{exist_set} $(iii)$, there exists
a $\mathbf{P}$-negligible set $\Omega_{\mathbf{f},\mathbf{t},0}(x)$
outside which the first two terms in the right-hand side of
(\ref{terme1bis}) converge to $0$. Besides, for any $\delta>0$,
\begin{eqnarray*}
\mathrm{E}_{0,x}\Big[\int_0^{\infty}ds\
\ind_{B}(\xi_s)|w_s^{*}(\xi_s)
- w_s^{(n)}(\xi_s)|\Big] & = & \int_0^{\infty}ds\ \int_{B}dz\ p_s(x,z)|w_s^{*}(z) - w_s^{(n)}(z)| \\
 & \leq & 2C_{\mathbf{f}}\delta +\frac{1}{(2\pi\delta)^{d/2}}\int_{\delta}^{\infty}ds\ \int_{B}dz\ |w_s^{*}(z)
- w_s^{(n)}(z)|,
\end{eqnarray*}
using the bound $p_s(x,z)\leq (2\pi \delta)^{-d/2}$ if $s\geq
\delta$. If in addition $\omega\in
\Big(\tilde{\Omega}_{\mathbf{f},\mathbf{t}} \cup
\hat{\Omega}_{\mathbf{f},\mathbf{t}}\Big)^c$, then by
(\ref{convergence1}) and dominated convergence (recall that $w_s^*$
and $w_s^{(n)}$ vanish for $s> t_p$),
$$\int_{\delta}^{\infty}ds\ \int_{B}dz\ |w_s^{*}(z) -
w_s^{(n)}(z)|\rightarrow 0$$ and so $\limsup
|w_0^{(n)}(x)-w_0^{*}(x)| \leq 2C_{\mathbf{f}}\delta$. Since
$\delta$ was arbitrary, it follows that $\lim
|w_0^{(n)}(x)-w_0^{*}(x)| =0$.

To summarize, for all $x\in B$ and $\omega\in
\Big(\tilde{\Omega}_{\mathbf{f},\mathbf{t}} \cup
\hat{\Omega}_{\mathbf{f},\mathbf{t}}\cup
\Omega_{\mathbf{f},\mathbf{t},0}(x)\Big)^c$ (of
$\mathbf{P}$-probability $1$), $$\lim_{n\rightarrow \infty}
|w_0^{(n)}(x)-w_0^{*}(x)| =0.$$

From the latter result, we can obtain the convergence of the
finite-dimensional distributions of $X^{(n)}$ towards the
corresponding ones for $X^*$. For all $x\in B$,
$\mathbf{P}\big[w_0^{(n)}(x)\rightarrow w_0^*(x)\big]=1$ so by
applying once again Fubini's theorem, we have
\begin{equation} \label{conv w}
\mathbf{P}\mathrm{-a.s.,}\quad \mu\mathrm{-a.e.,}\quad w_0^{(n)}(x)
\rightarrow w_0^*(x) \quad \mathrm{as}\ n\rightarrow \infty.
\end{equation}
Since the $w^{(n)}$ are bounded by $C_{\mathbf{f}}$, dominated
convergence and (\ref{conv w}) give $$\exp -\langle
\mu,w_0^{(n)}(\cdot)\rangle \rightarrow \exp -\langle
\mu,w_0^*(\cdot)\rangle.$$ Our construction from the historical
superprocess makes it obvious that $X^{(n)}$ is stochastically
bounded by $X^0$. It follows that the sequence of the distributions
of $\{(X^{(n)}_{t_1},\ldots,X^{(n)}_{t_p}), n\in \mathbb{N}\}$ is
relatively compact. Therefore, if we choose a countable set of
$p$-tuples $(f_1,\ldots,f_p)$ such that the corresponding family of
maps $(\mu_1,\ldots,\mu_p) \mapsto \exp - \sum_{i=1}^p \langle
\mu_i,f_i\rangle$ is convergence determining, we obtain that
$(X^{(n)}_{t_1},\ldots,X^{(n)}_{t_p})$ converges in distribution to
$(X^*_{t_1},\ldots,X^*_{t_p})$ on a set of $\mathbf{P}$-probability
$1$ (which a priori depends on $(t_1,\ldots,t_p)$). This completes
the proof of Proposition \ref{conv_fini_dim}.
\begin{flushright}$\Box$ \end{flushright}

\paragraph{Proof of Lemma \ref{borne essentielle}:} The quantity of interest vanishes
if $t>t_1$, and so we need only consider the case $t\leq t_1$. In
that case, \setlength\arraycolsep{1pt}
\begin{eqnarray}
\mathbf{E}&\bigg[&\Big(\mathrm{E}_{t,x}\Big[f(\xi_{t_1})\ind_{\{t_1
<
T\}}\Big(\ind_{\{t_1<T_{\e}\}}- e^{-k_d\int_t^{t_1}c(\xi_u)du}\Big)\Big]\Big)^2\bigg] \label{eq1 lemme 1}\\
& & =
\mathbf{E}\Big[\mathrm{E}_{t,x}\Big[\ind_{\{t_1<T\}}\ind_{\{t_1<T'\}}f(\xi_{t_1})f(\xi'_{t_1})\Big(
\ind_{\{t_1<T_{\e}\}}\ind_{\{t_1<T'_{\e}\}}-\ind_{\{t_1<T_{\e}\}}e^{-k_d\int_t^{t_1}c(\xi'_u)du}
\nonumber\\ & & \qquad \qquad \qquad -\
\ind_{\{t_1<T'_{\e}\}}e^{-k_d\int_t^{t_1}c(\xi_u)du}+
e^{-k_d\int_t^{t_1}\big(c(\xi_u)+c(\xi'_u)\big)du} \Big)\Big]\Big],
\nonumber
\end{eqnarray}
where $\xi'$ is another Brownian motion, independent of $\xi$, $T'$
and $T'_{\e}$ are defined in an obvious way and we have kept the
notation $\mathrm{P}_{t,x}$ for the probability measure under which
the two Brownian motions start from $x$ at time $t$. Recall that
$S_{\e}(s,t)$ denotes the Wiener sausage of radius $\e$ along the
time interval $[s,t]$ associated to $\xi$, and define $S'_{\e}(s,t)$
in a similar way. Then,
$$\mathbf{E}\big[\ind_{\{t_1<T_{\e}\}}\ind_{\{t_1<T'_{\e}\}}\big]=\mathbf{P}\big[\mathcal{P}^{\e}\cap
(S_{\e}(t,t_1)\cup S'_{\e}(t,t_1))=\emptyset\big]= \exp \Big\{ -
s_d(\e)\int_{S_{\e}(t,t_1)\cup S'_{\e}(t,t_1)}c(y)dy \Big\}$$ and
similarly
$$\mathbf{E}\big[\ind_{\{t_1<T_{\e}\}}\big]= \exp \Big\{-
s_d(\e)\int_{S_{\e}(t,t_1)}c(y)dy\Big\}. $$

Set $\tilde{t}_1=t_1-t$. By Fubini's theorem and a simple symmetry
argument, the quantity in (\ref{eq1 lemme 1}) is equal to
\setlength\arraycolsep{1pt}
\begin{eqnarray}
\mathrm{E}_{t,x}&\Big[&\ind_{\{t_1<T\}}\ind_{\{t_1<T'\}}f(\xi_{t_1})f(\xi'_{t_1})\Big\{
2e^{-k_d\int_t^{t_1}c(\xi'_u)du}\big(e^{-k_d\int_t^{t_1}c(\xi'_u)du}-e^{-s_d(\e)\int_{S_{\e}(t,t_1)}c(y)dy}\big)\nonumber
\\ & & +\ e^{-s_d(\e)\int_{S_{\e}(t,t_1)\cup
S'_{\e}(t,t_1)}c(y)dy}-e^{-k_d\int_t^{t_1}\big(c(\xi_u)+c(\xi'_u)\big)du}\Big\}\Big]
\nonumber \\
&\leq & 2\|f\|^2
\mathrm{E}_{t,x}\Big[\Big|e^{-k_d\int_t^{t_1}c(\xi'_u)du}-e^{-s_d(\e)\int_{S_{\e}(t,t_1)}c(y)dy}\Big|\Big]\nonumber
\\ & & +\|f\|^2
\mathrm{E}_{t,x}\Big[\Big|e^{-s_d(\e)\int_{S_{\e}(t,t_1)\cup
S'_{\e}(t,t_1)}c(y)dy}-e^{-k_d\int_t^{t_1}\big(c(\xi_u)+c(\xi'_u)\big)du}\Big|\Big]\nonumber \\
&\leq & 2\|f\|^2 \mathrm{E}_{0,x}\Big[\big|\ k_d
\int_0^{\tilde{t}_1}c(\xi_u)du -
s_d(\e)\int_{S_{\e}(0,\tilde{t}_1)}c(y)dy \big|\Big] + \|f\|^2
\mathrm{E}_{0,x}\Big[\big| s_d(\e)\int_{S_{\e}(0,\tilde{t}_1)}c(y)dy
\nonumber \\
& &+ s_d(\e)\int_{S'_{\e}(0,\tilde{t}_1)}c(y)dy -
s_d(\e)\int_{S_{\e}(0,\tilde{t}_1)\cap S'_{\e}(0,\tilde{t}_1)}c(y)dy
- k_d \int_0^{\tilde{t}_1}c(\xi_u)du - k_d
\int_0^{\tilde{t}_1}c(\xi'_u)du \big| \Big] \nonumber \\
&\leq & 4\|f\|^2 \mathrm{E}_{0,x}\Big[\big|\ k_d
\int_0^{\tilde{t}_1}c(\xi_u)du -
s_d(\e)\int_{S_{\e}(0,\tilde{t}_1)}c(y)dy \big|\Big] \nonumber \\& &
+ \|f\|^2\|c\|s_d(\e)
\mathrm{E}_{0,x}\big[\lambda\big(S_{\e}(0,\tilde{t}_1)\cap
S'_{\e}(0,\tilde{t}_1)\big)\big], \label{eq2 lemme 1}
\end{eqnarray}
where in the second inequality we used the bound
$|e^{-x}-e^{-y}|\leq |x-y|$ for $x,y\geq 0$.

On the one hand, by \cite{LEG1986} ($d=2,3$) and \cite{LEG1988}
p.1009-1010 ($d\geq 4)$, we have
\begin{eqnarray} |\log
\e|\ \mathrm{E}_{0,x}\big[\lambda\big(S_{\e}(0,\tilde{t}_1)\cap
S'_{\e}(0,\tilde{t}_1)\big)\big] &=& O(|\log \e|^{-1}) \quad
\mathrm{if\ }
d=2, \label{eq4 lemme 1}\\
\e^{-1}\mathrm{E}_{0,x}\big[\lambda\big(S_{\e}(0,\tilde{t}_1)\cap
S'_{\e}(0,\tilde{t}_1)\big)\big] &=& O(\e) \qquad \qquad \mathrm{if\ }d=3, \nonumber \\
\e^{2-d}\mathrm{E}_{0,x}\big[\lambda\big(S_{\e}(0,\tilde{t}_1)\cap
S'_{\e}(0,\tilde{t}_1)\big)\big] &=& O(\e^2 |\log \e|) \quad
\mathrm{if\ }d\geq 4. \nonumber
\end{eqnarray}
On the other hand,
\begin{equation}\label{eq3 lemme 1}\mathrm{E}_{0,x}\Big[\big|\ k_d
\int_0^{\tilde{t}_1}c(\xi_u)du -
s_d(\e)\int_{S_{\e}(0,\tilde{t}_1)}c(y)dy \big|\Big] \leq
\mathrm{E}_{0,x}\Big[\Big(\ k_d \int_0^{\tilde{t}_1}c(\xi_u)du -
s_d(\e)\int_{S_{\e}(0,\tilde{t}_1)}c(y)dy \Big)^2\Big]^{1/2}.
\end{equation}
Proposition \ref{prop section 2} ensures that the right-hand side of
(\ref{eq3 lemme 1}) is bounded by $K\|c\|\ |\log\e|^{-1}$ if $d=2$
and by $K\|c\|\ \e |\log\e|$ if $d\geq 3$. Together with (\ref{eq2
lemme 1}) and (\ref{eq4 lemme 1}), this completes the proof of Lemma
\ref{borne essentielle}.
\begin{flushright}$\Box$\end{flushright}

\section{Tightness of the sequence $X^{(n), B}$}
Let $C^2_+(\mathbb{R}^d)$ denote the set of all nonnegative
functions of class $C^2$ on $\mathbb{R}^d$. By Theorem II.4.1 in
\cite{PER1999}, the tightness of the sequence of the laws of the
superprocesses $X^{(n),B}$ will follow if we can prove that the
sequence of the laws of $\langle X^{(n),B},\phi\rangle$ is tight,
for every $\phi\in C^2_+(\mathbb{R}^d)$ with compact support. Note
that condition $(i)$ in Theorem II.4.1 of \cite{PER1999} holds
thanks to the domination $X^{(n),B}\leq X^0$. Recall that $X^0$ is
the usual super-Brownian motion without obstacles.

Let us first introduce the $\mathbf{P}$-negligible set
$\Theta\subset \Omega$ outside which the desired tightness will
hold.
\begin{defi}\textbf{(Good environments)} Let $\Theta$ be the union
over all choices of the integer $p\geq 1$ and of the rational
numbers $q_1,\ldots,q_p$ of the $\mathbf{P}$-negligible sets on
which the sequence $(X_{q_1}^{(n)},\ldots,X_{q_p}^{(n)})$ does not
converge in distribution to $(X_{q_1}^*,\ldots,X_{q_p}^*)$ as
$n\rightarrow \infty$. We call good environment any environment
which does not belong to $\Theta$.
\end{defi}

To simplify notation, we again write $X^{(n)}$ for $X^{(n),B}$ (as
in the last section, $B$ is fixed) and prove tightness only on the
time interval $[0,1]$. Let us fix $\phi \in C_+^2(\mathbb{R}^d)$
with compact support. The tightness of the sequence $\langle
X^{(n)},\phi\rangle$ is a consequence of the following lemma.
\begin{lemm}\label{uniform}If $\omega \notin \Theta$, then for every $\e>0$ there exists $k=k(\e)\geq 1$ and
$n_0=n_0(\mathbf{\omega},\e,k)$ such that
for all $n\geq n_0$,
$$\mathbb{P}_{\mu}\Big[\bigcup_{i=0}^{k-1}\Big\{\sup_{\frac{i}{k}\leq t\leq \frac{i+1}{k}}\big|\langle X_t^{(n)},
\phi\rangle -\langle X_{\frac{i}{k}}^{(n)},\phi \rangle\big|
>\e\Big\}\Big]<\e. $$
\end{lemm}

Lemma \ref{uniform} easily implies that the sequence $\langle
X^{(n)},\phi\rangle$ is tight. Indeed, let us fix a good environment
and $\eta>0$. By Lemma \ref{uniform}, there exist $k(\eta)$ and
$n_0(\omega,\eta,k)$ such that for all $n\geq n_0$,
\begin{equation}\label{tightness ineq1}
\mathbb{P}_{\mu}\Big[\bigcup_{i=0}^{k-1}\Big\{\sup_{\frac{i}{k}\leq
t\leq \frac{i+1}{k}}\big|\langle X_t^{(n)},\phi\rangle -\langle
X_{\frac{i}{k}}^{(n)},\phi \rangle\big|
>\frac{\eta}{3}\Big\}\Big]<\eta. \end{equation}
If $n\geq n_0$ is fixed, on the complement of the event considered
in (\ref{tightness ineq1}), we have for every $s,t\in [0,1]$
$$|t-s|\leq \frac{1}{k}\quad \Rightarrow \quad \big|\langle
X_t^{(n)},\phi\rangle -\langle X_s^{(n)},\phi \rangle \big|\leq \eta
$$ and therefore $w(\langle X^{(n)},\phi\rangle,\frac{1}{k},1)\leq
\eta$, using the notation of Ethier and Kurtz \cite{EK1986} for the
modulus of continuity of the process $\langle X^{(n)},\phi\rangle$.
Thus, for all $n\geq n_0$,
$$
\mathbb{P}_{\mu}\big[w(\langle
X^{(n)},\phi\rangle,\frac{1}{k},1)\leq \eta \big] \geq
\mathbb{P}_{\mu}\Big[\Big(\bigcup_{i=0}^{k-1}\Big\{\sup_{\frac{i}{k}\leq
t\leq \frac{i+1}{k}}\big|\langle X_t^{(n)},\phi\rangle -\langle
X_{\frac{i}{k}}^{(n)},\phi \rangle\big|
>\frac{\eta}{3}\Big\}\Big)^c\Big] \geq 1-\eta. $$
In addition, $\phi$ is bounded so the first condition of Theorem
3.7.2 in \cite{EK1986} is trivially fulfilled, hence Corollary 3.7.4
of \cite{EK1986} implies that for any good environment, the sequence
of the laws of $\langle X^{(n)},\phi\rangle$ under
$\mathbb{P}_{\mu}$ is tight.

Let us now turn to the proof of Lemma \ref{uniform}.
\paragraph{Proof of Lemma \ref{uniform}:} We fix a good
environment. Let $\e>0$. The process $(\langle X^*_t,\phi
\rangle)_{t\geq 0}$ is continuous, therefore there exists $k_0(\e)$
such that for all $k\geq k_0$,
\begin{equation}\label{cont X*}
\mathbb{P}_{\mu}\Big[\sup_{0\leq i\leq k-1}\Big|\langle
X^*_{\frac{i+1}{k}},\phi\rangle - \langle
X^*_{\frac{i}{k}},\phi\rangle\Big|\geq \frac{\e}{2}\Big] <
\frac{\e}{3}. \end{equation} There exists $K=K(\e)\geq 1$ such that
\begin{equation}\label{cont masse de X}
\mathbb{P}_{\mu}\Big[\sup_{0\leq t\leq 1}\langle X_t^0,1\rangle \geq
K\Big] < \frac{\e}{3}.
\end{equation}
By a trivial domination argument, the bound (\ref{cont masse de X})
remains valid if we replace $X^0$ by $X^{(n)}$ (in fact for any
environment). In the following, we fix the constant $K\geq 1$ such
that (\ref{cont masse de X}) holds.

We now have the following result :
\begin{lemm}\label{cont X}There exists a constant $C=C(\phi,K)$ such that for every integer $k\geq 1$ and
every measure $\gamma\in \mathcal{M}_f(\mathbb{R}^d)$ satisfying
$\langle \gamma, 1\rangle \leq K$,
$$\mathbb{P}_{\gamma}\Big[\sup_{0\leq s\leq \frac{1}{k}}\Big|\langle X_s^0,\phi \rangle -\langle
X_0^0,\phi\rangle\Big|>\frac{\e}{2}\Big] \leq \frac{C}{k^2}\ .$$
\end{lemm}

The proof of Lemma \ref{cont X} is deferred to the end of the
section. Let us define
$$A_n=\big\{\sup_{0\leq t\leq 1}\langle X_t^{(n)},1\rangle
\geq K \big\}.$$ Then, \setlength\arraycolsep{1pt}
\begin{eqnarray*}
&\mathbb{P}_{\mu}&\Big[\bigcup_{i=0}^{k-1}\Big\{\sup_{\frac{i}{k}\leq
t\leq \frac{i+1}{k}}\big|\langle X_t^{(n)},\phi\rangle -\langle
X_{\frac{i}{k}}^{(n)},\phi \rangle\big|
>\e\Big\}\Big] \\ &\leq& \mathbb{P}_{\mu}[A_n] + \mathbb{P}_{\mu}\Big[\sup_{0\leq i\leq k-1}\Big|\langle
X^{(n)}_{\frac{i+1}{k}},\phi\rangle - \langle
X^{(n)}_{\frac{i}{k}},\phi\rangle\Big|> \frac{\e}{2}\Big] \\
& +& \mathbb{P}_{\mu}\Big[A_n^c \cap \Big\{\sup_{0\leq i\leq
k-1}\Big\{\sup_{\frac{i}{k}\leq t \leq \frac{i+1}{k}}\Big(\langle
X^{(n)}_t,\phi\rangle - \langle
X^{(n)}_{\frac{i}{k}},\phi\rangle\Big)\Big\}>\e\Big\}\Big] \\
& +& \mathbb{P}_{\mu}\Big[A_n^c \cap \Big\{\sup_{0\leq i\leq
k-1}\Big\{\sup_{\frac{i}{k}\leq t \leq \frac{i+1}{k}}\Big(\langle
X^{(n)}_{\frac{i}{k}},\phi\rangle - \langle
X^{(n)}_t,\phi\rangle\Big)\Big\}>\e\Big\} \cap \Big\{\sup_{0\leq
i\leq k-1}\Big|\langle X^{(n)}_{\frac{i+1}{k}},\phi\rangle - \langle
X^{(n)}_{\frac{i}{k}},\phi\rangle\Big|\leq \frac{\e}{2}\Big\}\Big] \\
&=& a_n+b_n+c_n+d_n
\end{eqnarray*}
From (\ref{cont masse de X}), we have $$a_n < \frac{\e}{3}.$$
Moreover, from the definition of a good environment, if $k\geq k_0$,
$$\limsup_{n\rightarrow \infty} b_n \leq \mathbb{P}_{\mu}\Big[\sup_{0\leq i\leq
k-1}\Big|\langle X^*_{\frac{i+1}{k}},\phi\rangle - \langle
X^*_{\frac{i}{k}},\phi\rangle\Big|\geq \frac{\e}{2}\Big] <
\frac{\e}{3},$$ by (\ref{cont X*}). Thus if $k\geq k_0$, there
exists $n_0(\omega,\e, k)$ such that for all $n\geq n_0$,
$b_n(k)\leq \frac{\e}{3}$. Then,
\begin{eqnarray*} c_n &\leq& \sum_{i=0}^{k-1}\mathbb{P}_{\mu}\Big[\langle X^{(n)}_{\frac{i}{k}},1\rangle \leq K;
\ \sup_{\frac{i}{k}\leq t \leq \frac{i+1}{k}}\Big(\langle
X^{(n)}_t,\phi\rangle - \langle
X^{(n)}_{\frac{i}{k}},\phi\rangle\Big)>\e \Big] \\
&=& \sum_{i=0}^{k-1}\mathbb{E}_{\mu}\Big[\ind_{\big\{\langle
X^{(n)}_{\frac{i}{k}},1\rangle
  \leq K \big\} }\ \mathbb{P}_{X_{\frac{i}{k}}^{(n)}}\Big[\sup_{0\leq t\leq \frac{1}{k}}\big(\langle
X^{(n)}_t,\phi\rangle - \langle X^{(n)}_0,\phi\rangle\big)>\e
\Big]\Big]
\end{eqnarray*}
The last equality is obtained by applying the Markov property to
$X^{(n)}$ at time $\frac{i}{k}$. By a domination argument, we have
for all $\gamma\in \mathcal{M}_f(\mathbb{R}^d)$ such that $\langle
\gamma,1\rangle \leq K$,
$$\mathbb{P}_{\gamma}\Big[\sup_{0\leq t\leq \frac{1}{k}}\big(\langle
X^{(n)}_t,\phi\rangle - \langle X^{(n)}_0,\phi\rangle\big)>\e \Big]
\leq \mathbb{P}_{\gamma}\Big[\sup_{0\leq t\leq
\frac{1}{k}}\big(\langle X^0_t,\phi\rangle - \langle
X^0_0,\phi\rangle\big)>\e \Big]\leq \frac{C}{k^2}$$ by Lemma
\ref{cont X}. It follows that
$$c_n\leq k.\frac{C}{k^2} =
\frac{C}{k} < \frac{\e}{6}$$ if $k\geq k_1(\e)$. Finally,
\begin{eqnarray*}
d_n &\leq &
\sum_{i=0}^{k-1}\mathbb{P}_{\mu}\Big[\sup_{\frac{i}{k}\leq t \leq
\frac{i+1}{k}}\langle X^{(n)}_t,1\rangle \leq K;\
\sup_{\frac{i}{k}\leq t \leq \frac{i+1}{k}}\Big(\langle
X^{(n)}_{\frac{i}{k}},\phi\rangle - \langle
X^{(n)}_t,\phi\rangle\Big)>\e; \\
& & \qquad \qquad \sup_{0\leq i\leq k-1}\Big|\langle
X^{(n)}_{\frac{i+1}{k}},\phi\rangle - \langle
X^{(n)}_{\frac{i}{k}},\phi\rangle\Big|\leq \frac{\e}{2}\Big].
\end{eqnarray*}
We fix $i\in \{0,\ldots, k-1\}$ and consider the stopping time
$$T_i:=\inf\Big\{t\geq \frac{i}{k}:\ \langle X_t^{(n)},\phi\rangle \leq
\langle X_{\frac{i}{k}}^{(n)},\phi\rangle - \e \Big\}.$$ Then, the
$i$-th term of the previous sum is bounded by
\setlength\arraycolsep{1pt}
\begin{eqnarray*}\mathbb{P}_{\mu}&&\Big[\ T_i \leq \frac{i+1}{k},\ \langle
  X^{(n)}_{T_i},1\rangle \leq K, \sup_{T_i\leq t \leq T_i+\frac{1}{k}}\Big(\langle
X^{(n)}_t,\phi\rangle - \langle
X^{(n)}_{T_i},\phi\rangle\Big)\geq \frac{\e}{2} \Big] \\
& & = \mathbb{E}_{\mu}\Big[\ind_{\big\{T_i\leq\frac{i+1}{k},\
\langle X^{(n)}_{T_i},1\rangle \leq K\big\}}\
\mathbb{P}_{X^{(n)}_{T_i}}\Big[\sup_{0\leq t\leq
\frac{1}{k}}\big(\langle X^{(n)}_t,\phi\rangle - \langle
X^{(n)}_0,\phi\rangle\big)\geq \frac{\e}{2} \Big]\Big]
\end{eqnarray*}
by the strong Markov property at time $T_i$. Using Lemma \ref{cont
X} once again, we see that this quantity is bounded by $C/k^2$,
hence for $k\geq k_1(\e)$,$$d_n \leq \frac{C}{k} \leq
\frac{\e}{6}.$$ Combining the preceding estimates, we obtain that
for $k = k_0(\e)\vee k_1(\e)$, and every $n\geq n_0(\omega,\e,k)$,
$$\mathbb{P}_{\mu}\Big[\bigcup_{i=0}^{k-1}\Big\{\sup_{\frac{i}{k}\leq
t\leq \frac{i+1}{k}}\big|\langle X_t^{(n)},\phi\rangle -\langle
X_{\frac{i}{k}}^{(n)},\phi \rangle\big|
>\e\Big\}\Big]<\e. $$
This completes the proof of Lemma \ref{uniform}.
\begin{flushright} $\Box$ \end{flushright}

\paragraph{Proof of Lemma \ref{cont X}:} Let $\gamma\in
\mathcal{M}_f(\mathbb{R}^d)$ be such that $|\gamma|:= \langle
\gamma,1\rangle \leq K$. Recall that the process $(\langle
X_t^0,1\rangle)_{t\geq 0}$ is a martingale. From the maximal
inequality applied to the nonnegative submartingale $(\langle
X_t^0,1\rangle - |\gamma|)^4$,
$$\mathbb{P}_{\gamma}\Big[\sup_{0\leq t \leq \frac{1}{k}}\langle X_t^0,1\rangle > 2K \Big]\leq \frac{1}{(2K-|\gamma|)^4}
\mathbb{E}_{\gamma}\Big[\big(\langle X^0_{\frac{1}{k}},1\rangle -
|\gamma|\big)^4 \Big]. $$ We claim that
\begin{equation}\label{moment csbp}\mathbb{E}_{\gamma}\Big[\big(\langle X^0_{\frac{1}{k}},1\rangle
- |\gamma|\big)^4 \Big] =
\frac{24}{k^3}|\gamma|+\frac{12}{k^2}|\gamma|^2.\end{equation} To
prove this claim, recall that $Y_t=\langle X_t^0,1\rangle$ is a
Feller diffusion, whose semigroup Laplace transform is given by
$$\mathbb{E}\big[\exp -\lambda Y_t |\ Y_0=y\big]=
\exp\Big(-\frac{\lambda y}{1+\lambda t}\Big),$$ for $\lambda\geq 0$.
Thus,
$$
\mathbb{E}\big[\exp -\lambda (Y_t-y) |\ Y_0=y\big]=
\exp\Big(\frac{\lambda^2 t y}{1+\lambda t}\Big)= 1+\lambda^2ty
-\lambda^3 t^2y+ \lambda^4t^3y
+\frac{\lambda^4t^2y^2}{2}+o(\lambda^4),
$$ as $\lambda \rightarrow 0$. From this expansion of the Laplace
transform, we derive that $$\mathbb{E}\Big[\big(Y_t-y)^4|Y_0=y \Big]
= 24t^3y+12t^2y^2,$$ which proves our claim (\ref{moment csbp}). It
follows that
$$\mathbb{P}_{\gamma}\Big[\sup_{0\leq t \leq \frac{1}{k}}\langle
X^0_t,1\rangle > 2K \Big]\leq
\frac{12|\gamma|(|\gamma|+2)}{(2K-|\gamma|)^4\ k^2}.$$

Let us denote by $A_{K,k}$ the event $\big\{\sup_{0\leq t \leq
\frac{1}{k}}\langle X_t^0,1\rangle > 2K\big\}$ and by $B_{K,k}$ the
event $\big\{\sup_{0\leq t\leq \frac{1}{k}}\big|\langle
X_t^0,\phi\rangle - \langle
X_0^0,\phi\rangle\big|>\frac{\e}{2}\big\}$. Then,
\begin{equation}\label{borne BKk}
\mathbb{P}_{\gamma}[B_{K,k}] \leq \mathbb{P}_{\gamma}[A_{K,k}] +
\mathbb{P}_{\gamma}[A_{K,k}^c  \cap B_{K,k}] \leq
\frac{\mathrm{c_0}}{k^2}+\mathbb{P}_{\gamma}[A_{K,k}^c  \cap
B_{K,k}],
\end{equation}
where $c_0$ is a constant depending on $K$.

In addition, $M_t:= \langle X_t^0,\phi\rangle-\langle
X_0^0,\phi\rangle - \int_0^tdr\ \langle
X_r^0,\frac{1}{2}\Delta\phi\rangle$ is a continuous martingale with
quadratic variation $2\int_0^tdr\ \langle X_r^0,\phi^2\rangle$. By
the Dubins-Schwarz theorem (see Theorem V.1.7 in \cite{RY1999}),
there exists a standard one-dimensional Brownian motion $W$ such
that $M_t= W_{\langle M\rangle_t}$ for all $t\geq 0$ a.s. On the
event $A_{K,k}^c$, we have
$$\Big|\int_0^t dr\ \langle X_r^0,\frac{1}{2}\Delta \phi\rangle \Big| \leq t\|\Delta\phi\|K \leq \frac{c_1}{k}
\quad\mathrm{if}\ t\in[0,k^{-1}]$$ and $$\langle M\rangle_t \leq
\frac{4\|\phi\|^2K}{k}=\frac{c_2}{k}\quad\mathrm{if}\ t\in
[0,k^{-1}],$$ where $c_1$ and $c_2$ are constants depending on
$\phi$ and on $K$. Choose $k_0$ such that $c_1 k^{-1}<\frac{1}{4}\e$
for every $k\geq k_0$. Then for all $k\geq k_0$,
\begin{eqnarray*}
\mathbb{P}_{\gamma}[A_{K,k}^c  \cap B_{K,k}] &=&
\mathbb{P}_{\gamma}\Big[A_{K,k}^c  \cap \big\{\sup_{0\leq t\leq 1/k}
\big|M_t+\int_0^t dr\ \langle X_r^0,\frac{1}{2}\Delta\phi\rangle \big|>\frac{\e}{2}\big\}\Big] \\
&\leq& \mathbb{P}_{\gamma}\Big[A_{K,k}^c  \cap \big\{\sup_{0\leq t\leq 1/k}|M_t|>\frac{\e}{4}\big\}\Big] \\
&\leq& \mathrm{P}\Big[\sup_{0\leq t\leq (c_2/k)}|W_t|>\frac{\e}{4}\Big] \\
&\leq &\frac{c_3}{k^2}\ ,
\end{eqnarray*}
where $c_3$ is a constant depending on $\phi$, $K$ and $\e$.
Together with (\ref{borne BKk}), this completes the proof of Lemma
\ref{cont X}.
\begin{flushright}$\Box$ \end{flushright}

\section{Proof of Theorem \ref{main theorem} and Corollary \ref{corollaire}}
The proof of Theorem \ref{main theorem} in the case when $D$ is
bounded is easy from the results of the previous sections. Let us
take $B=D$ and let $\mathcal{E}$ denote the union of the
$\mathbf{P}$-negligible set on which there exist rational numbers
$t_1,\ldots,t_p$ such that $(X_{t_1}^{(n)},\ldots,X_{t_p}^{(n)})$
does not converge to $(X_{t_1}^*,\ldots,X_{t_p}^*)$ and of the
$\mathbf{P}$-negligible set on which the sequence $X^{(n),B}$ is not
tight. The set $\mathcal{E}$ is also $\mathbf{P}$-negligible and on
$\mathcal{E}^c$, Theorem 3.7.8 of \cite{EK1986} allows us to
conclude that $X^{(n)} \stackrel{(d)}{\rightarrow} X^*$ when
$n\rightarrow \infty$.

We can now use the previous result to complete the proof of Theorem
\ref{main theorem} when $D$ is a domain of $\mathbb{R}^d$ which is
not necessarily bounded.
\paragraph{Proof of Theorem \ref{main theorem} for a general domain $D$:}
Let $\mu$ be a finite measure on $D$ and suppose first that the
support of $\mu$ is bounded. Under $\mathbb{P}_{\mu}$, the
superprocesses $X^{(n)}$ are stochastically dominated by the
superprocess $X^0$, whose range
$$\mathcal{R}(X^0)=\overline{\bigcup_{t \geq\ 0}\mathrm{supp}X_t^0}$$ is almost surely compact
since its initial value has compact support. Consequently, for every
$\e>0$, there exists a bounded open subset $B$ of $D$ containing the
support of $\mu$ such that, for every environment and every $n\geq
1$,
$$\mathbb{P}_{\mu}[\mathcal{R}(X^{(n),D})\subset B]\geq 1-\e$$
and $$\mathbb{P}_{\mu}[\mathcal{R}(X^{*,D})\subset B]\geq 1-\e.$$
From these inequalities, we can deduce that
\begin{eqnarray}
d(\mathbb{P}_{\mu}^{(n),D},\mathbb{P}_{\mu}^{(n),B}) &\leq& 2\e,\qquad n\geq 1; \label{distance_n} \\
d(\mathbb{P}_{\mu}^{*,D},\mathbb{P}_{\mu}^{*, B}) &\leq& 2\e,
\label{distance_*}
\end{eqnarray}
where $d$ is the Prohorov metric on
$\mathcal{M}_1(D_{\mathcal{M}_f(D)}[0,\infty))$. By the results of
the last two sections, with $\mathbf{P}$-probability $1$ there
exists an integer $n_0(\omega)$ such that for all $n\geq n_0$,
$$d(\mathbb{P}_{\mu}^{(n), B},\mathbb{P}_{\mu}^{*, B})\leq \e.$$
Together with (\ref{distance_n}) and (\ref{distance_*}), this yields
$$d(\mathbb{P}_{\mu}^{(n),D},\mathbb{P}_{\mu}^{*,D})\leq 5\e \qquad \qquad \mathrm{for\ all\ }n \geq n_0(\omega),$$
hence we can conclude that $\mathbb{P}_{\mu}^{(n),D}$ converges
towards $\mathbb{P}_{\mu}^{*,D}$ on a set of
$\mathbf{P}$-probability $1$.

Finally, if the support of $\mu$ is unbounded, we can replace $\mu$
by the measure $\tilde{\mu}$ defined as the restriction of $\mu$ to
a large ball centered at the origin. Using once again the domination
of $X^{(n),D}$ (for all $n\geq 1$) and of $X^{*,D}$ by $X^0$, the
law of $X^{(n),D}$ under $\mathbb{P}_{\mu}$ can be approximated
uniformly in $n$ by the law of $X^{(n),D}$ under
$\mathbb{P}_{\tilde{\mu}}$, and similarly for $X^{*,D}$. The desired
result then follows from the bounded support case. We leave the
details to the reader.
\begin{flushright}$\Box$ \end{flushright}
We end this section with the proof of Corollary \ref{corollaire}.

\paragraph{Proof of Corollary \ref{corollaire}:}
Let us argue by contradiction and suppose that there exist $\delta
>0$ and a sequence $\{\e_k,k\in \mathbb{N}\}$ decreasing to zero
such that for all $k\geq 1$,
\begin{equation}\label{contradict}\mathbf{P}\left[d\big(\mathbb{P}^{\e_k,D}_{\mu},\mathbb{P}^{*,D}_{\mu}\big)>
\delta \right] \geq \delta. \end{equation}
By extracting a
subsequence, we can always choose $\e_k$ such that
$$\sum_{k=1}^{\infty}|\log \e_k|^{-1} <\infty \mathrm{\;\; if\;
}d=2,$$ or $$\sum_{k=1}^{\infty}\e_k\ |\log \e_k| < \infty
\mathrm{\;\; if\; }d\geq 3.$$ But the latter condition is the only
requirement for the sequence of superprocesses $X^{\e_k,D}$ to
converge in distribution to $X^{*,D}$ with $\mathbf{P}$-probability
$1$, yielding a contradiction with (\ref{contradict}).
\begin{flushright}$\Box$ \end{flushright}

\appendix
\section{Appendix: Proof of (\ref{borne v,2})}
The bound (\ref{borne v,2}) is a consequence of the following lemma.
\begin{lemm} \label{lemma appendix}
There exists a function $\varphi:\mathbb{R}^2\rightarrow [0,\infty]$
such that $\int_{\mathbb{R}^2}dy\ \varphi(y)<\infty$ and for every
$y\in \mathbb{R}^2$ and $\e\in (0,\frac{1}{2})$,
$$\Big|\mathrm{P}_0\big[y\in S_{\e}(0,1)\big]-\frac{\pi}{|\log\e|}\int_0^1ds\
p_s(y)\Big|\leq \frac{\varphi(y)}{|\log\e|^2},$$ where $p_s(y)=(2\pi
s)^{-1}\exp\big\{-|y|^2/(2s)\big\}$ is the Brownian transition
density.
\end{lemm}
\begin{rema}The convergence of $|\log\e| \mathrm{P}_0[y\in S_{\e}(0,1)]$ towards $\pi \int_0^1 ds\ p_s(y)$
as $\e$ tends to $0$ was first obtained by Spitzer \cite{SPI1958}.
See also \cite{LEG1986} and \cite{SZN1987} for related results.
\end{rema}

Before proving Lemma \ref{lemma appendix}, let us use it to derive
(\ref{borne v,2}). If $c$ is a bounded nonnegative measurable
function on $\mathbb{R}^2$ such that $\|c\|\leq 1$, then for every
$\e \in (0,\frac{1}{2}]$ and $z\in \mathbb{R}^2$,
\begin{eqnarray*}\Big|\mathrm{E}_z\Big[|\log\e|\int_{S_{\e}(0,1)}c(y)dy-\pi\int_0^1c(\xi_s)ds\Big]\Big|
& =&\Big|\int dy\ c(z+y)\Big(|\log\e|\mathrm{P}_0\big[y\in
S_{\e}(0,1) \big]-\pi \int_0^1ds\ p_s(y)\Big)\Big|\\
&\leq & |\log\e|^{-1}\int dy\ \varphi(y),
\end{eqnarray*}
which is the desired result.

Let us hence establish Lemma \ref{lemma appendix}. The following
proof is due to J.-F. Le Gall (personal communication).

\medskip
\noindent\textbf{Proof of Lemma \ref{lemma appendix}: }If $|y|\leq
10\e$, simple estimates show that the bound of Lemma \ref{lemma
appendix} holds with $\varphi(y)=C\big(\big(\log |y|\big)^2+1\big)$
for a suitable constant $C$. So we assume that $|y|>10\e$. We put
$$\tau_{\e}(y)=\inf\big\{t\geq 0:|\xi_t-y|\leq \e\big\}$$
in such a way that $\{y\in S_{\e}(0,1)\}=\{\tau_{\e}(y)\leq 1\}$.
Let $a_{\e}$ be an arbitrary point of the circle of radius $\e$
centered at the origin, and
$$f(\e)= \mathrm{E}_{a_{\e}}\bigg[\int_0^1ds\ \ind_{\{|\xi_s|\leq
\e\}}\bigg].$$ A straightforward calculation gives
\begin{equation}\label{eq appendix}f(\e)=\e^2|\log\e|+O(\e^2)\end{equation} as $\e\rightarrow 0$.

\noindent\emph{Lower bound.} An application of the strong Markov
property at time $\tau_{\e}(y)$ shows for every $u\in (0,1]$, that
$$\mathrm{E}_0\bigg[\int_0^uds\ \ind_{\{|\xi_s-y|\leq \e\}}\bigg]\leq \mathrm{P}_0\big[\tau_{\e}(y)\leq u\big]f(\e).$$
On the other hand,
$$\mathrm{E}_0\bigg[\int_0^uds\ \ind_{\{|\xi_s-y|\leq \e\}}\bigg]= \int_0^u ds\ \int_{|z-y|\leq \e}dz \ p_s(z),$$
and thus \setlength\arraycolsep{1pt}
\begin{eqnarray}
\bigg|\mathrm{E}_0\bigg[\int_0^u&ds&\ \ind_{\{|\xi_s-y|\leq
\e\}}\bigg]-\pi \e^2\int_0^u ds\ p_s(y)\bigg|\nonumber\\
&\leq & \int_0^u \frac{ds}{2\pi s}\int_{|z-y|\leq \e}dz\
\bigg|\exp\Big\{-\frac{|z|^2}{2s}\Big\}
-\exp\Big\{-\frac{|y|^2}{2s}\Big\}\bigg|\nonumber\\
&\leq & \int_0^u \frac{ds}{2\pi s}\int_{|z-y|\leq \e}dz\ \exp\Big\{-\frac{|y|^2}{4s}\Big\}
\bigg|\frac{|z|^2-|y|^2}{2s}\bigg|\nonumber\\
&\leq & \frac{\e^3}{2}|y|\int_0^u\frac{ds}{s^2}\ \exp\Big\{-\frac{|y|^2}{4s}\Big\}\nonumber\\
& \leq & \e^3\ \Psi_1(y), \label{eq1 lower bound}
\end{eqnarray}
where the function
$$\Psi_1(y)=|y|\int_0^1\frac{ds}{s^2}\ \exp\Big\{-\frac{|y|^2}{4s}\Big\}$$ is
easily seen to be integrable over $\mathbb{R}^2$.

By combining the preceding estimates, we arrive at
\begin{equation}\label{eq2 lower bound}
\mathrm{P}_0\big[\tau_{\e}(y)\leq u\big]\geq \frac{\pi
\e^2}{f(\e)}\int_0^uds\ p_s(y) -\frac{\e^3}{f(\e)}\Psi_1(y)
\end{equation} and using (\ref{eq appendix}) it readily follows that
$$\mathrm{P}_0\big[y\in S_{\e}(0,1)\big]-\frac{\pi}{|\log\e|}\int_0^1 ds\ p_s(y)\geq \frac{\varphi_1(y)}{|\log\e|^2}$$
with a nonnegative function $\varphi_1$ such that $\int dy\
\varphi_1(y)<\infty$.

\medskip
\noindent\emph{Upper bound.} This part is a little more delicate. We
rely on the same idea of applying the strong Markov property at time
$\tau_{\e}(y)$, but we need to be more careful in our estimates. For
every $v>0$, we have \setlength\arraycolsep{1pt}
\begin{eqnarray*}
\mathrm{E}_0\bigg[\int_0^{1+v}ds\ \ind_{\{|\xi_s-y|\leq \e\}}\bigg]
 &=& \mathrm{E}_0\bigg[\ind_{\{\tau_{\e}(y)\leq
1+v\}}\mathrm{E}_{\xi_{\tau_{\e}(y)}}
\bigg[\int_0^s dr\ \ind_{\{|\xi_r-y|\leq \e\}}\bigg]_{s=1+v-\tau_{\e}(y)}\bigg]\\
&=&\mathrm{E}_0\bigg[\ind_{\{\tau_{\e}(y)\leq
1+v\}}\int_0^{1+v-\tau_{\e}(y)}dr\
\mathrm{P}_{a_{\e}}\big[|\xi_r|\leq \e\big]\bigg],
\end{eqnarray*}
where $a_{\e}$ is as previously a fixed point of the circle of
radius $\e$ centered at the origin. We can rewrite the previous
expression as
$$\mathrm{E}_0\bigg[\int_0^{1+v}dr\ \ind_{\{\tau_{\e}(y)\leq 1+v-r\}}\mathrm{P}_{a_{\e}}\big[|\xi_r|\leq \e\big]\bigg]=
\int_0^{1+v}dr\ \mathrm{P}_0\big[\tau_{\e}(y)\leq
1+v-r\big]\mathrm{P}_{a_{\e}}\big[|\xi_r|\leq \e\big].$$

We apply this calculation with $v=v_{\e}=|\log\e|^{-1}$. For $r\in
[0,v_{\e}]$, $\mathrm{P}_0\big[\tau_{\e}(y)\leq 1+v_{\e}-r\big]$ is
bounded from below by $\mathrm{P}_0\big[\tau_{\e}(y)\leq 1\big]$,
and thus \setlength\arraycolsep{1pt}
\begin{eqnarray*}\mathrm{P}_0&\big[&\tau_{\e}(y)\leq 1\big]\int_0^{v_{\e}}dr\
\mathrm{P}_{a_{\e}}\big[|\xi_r|\leq \e\big]\\
&\leq &\mathrm{E}_0\Big[\int_0^{1+v_{\e}}ds\ \ind_{\{|\xi_s-y|\leq
\e\}}\Big]-\int_{v_{\e}}^{1+v_{\e}}dr\
\mathrm{P}_{a_{\e}}\big[|\xi_r|\leq
\e\big]\mathrm{P}_0\big[\tau_{\e}(y)\leq 1+v_{\e}-r\big].
\end{eqnarray*}

From the bound (\ref{eq1 lower bound}), we have
$$\mathrm{E}_0\bigg[\int_0^{1+v_{\e}}ds\ \ind_{\{|\xi_s-y|\leq \e\}}\bigg]
\leq \pi \e^2\int_0^1 ds\ p_s(y)+\e^3\Psi_1(y)+v_{\e}\e^2
e^{-|y|^2/10}.$$ On the other hand, by (\ref{eq2 lower bound}),
\setlength\arraycolsep{1pt}
\begin{eqnarray*}\int_{v_{\e}}^{1+v_{\e}}&dr& \mathrm{P}_{a_{\e}}\big[|\xi_r|\leq \e\big]
\mathrm{P}_0\big[\tau_{\e}(y)\leq 1+v_{\e}-r\big]\\
&\geq & \int_{v_{\e}}^{1+v_{\e}}dr\ \mathrm{P}_{a_{\e}}\big[|\xi_r|\leq \e\big]
\frac{\e^2}{f(\e)}\Big(\pi\int_0^{1+v_{\e}-r}ds\ p_s(y)-\e\Psi_1(y)\Big)\\
&= &\frac{\pi \e^2}{f(\e)}\Big(\int_{v_{\e}}^{1+v_{\e}}dr\
\mathrm{P}_{a_{\e}}\big[|\xi_r|\leq \e\big]\Big)\Big(\int_0^1ds\
p_s(y)-\frac{\e}{\pi} \Psi_1(y)\Big)\\
& & \quad -\frac{\pi\e^2}{f(\e)}\Big(\int_{v_{\e}}^{1+v_{\e}}dr\
\mathrm{P}_{a_{\e}}\big[|\xi_r|\leq \e\big]\int_{1+v_{\e}-r}^1ds\
p_s(y)\Big).
\end{eqnarray*}
Straightforward estimates give
$$\int_{v_{\e}}^{1+v_{\e}}dr\ \mathrm{P}_{a_{\e}}\big[|\xi_r|\leq
\e\big]=\e^2\Big(\frac{1}{2}\log|\log\e|+O(1)\Big)$$ and
$$\int_{v_{\e}}^{1+v_{\e}}dr\ \mathrm{P}_{a_{\e}}\big[|\xi_r|\leq
\e\big]\int_{1+v_{\e}-r}^1ds \ p_s(y)\leq \e^2\Psi_2(y),$$ where
$$\Psi_2(y)=\int_0^1ds\ \log\Big(\frac{1}{1-s}\Big)p_s(y)$$ is
integrable over $\mathbb{R}^2$. Summarizing, we have
\setlength\arraycolsep{1pt}
\begin{eqnarray*}
\mathrm{P}_0\big[\tau_{\e}(y)\leq 1\big]\int_0^{v_{\e}}dr\
\mathrm{P}_{a_{\e}}\big[|\xi_r|\leq \e\big] &\leq &\Big(\pi
\e^2\int_0^1ds\
p_s(y)\Big)\times\Big(1-\frac{\big(\frac{1}{2}\log|\log\e|-K\big)
\e^2}{f(\e)}\Big)\\
& & +
\Big(\e^3+O\Big(\e^3\frac{\log|\log\e|}{|\log\e|}\Big)\Big)\Psi_1(y)+v_{\e}\e^2e^{-|y|^2/10}+\frac{\pi
\e^4}{f(\e)}\Psi_2(y).
\end{eqnarray*}
Finally, it is easy to verify that
$$g(\e)\equiv \int_0^{v_{\e}}dr\ \mathrm{P}_{a_{\e}}\big[|\xi_r|\leq
\e\big]\geq \e^2\big(|\log\e| -\frac{1}{2}\log|\log\e|-K'\big),$$
and so by dividing the two sides of the previous inequality by
$g(\e)$, we obtain
$$\mathrm{P}_0\big[\tau_{\e}(y)\leq 1\big]\leq
\frac{\pi}{|\log\e|}\int_0^1ds\
p_s(y)+\frac{\varphi_2(y)}{|\log\e|^2},$$ with a function
$\varphi_2$ such that $\int \varphi_2(y)dy<\infty$. This completes
the proof of Lemma \ref{lemma appendix}.
\begin{flushright}
$\Box$
\end{flushright}

\bigskip
\noindent \textbf{Acknowledgements.} The author would like to thank
her supervisor Jean-Fran\c cois Le Gall for many helpful discussions
about this work and detailed comments on earlier drafts of this
paper, and the referee for valuable suggestions to improve the
presentation of the paper.

\end{document}